\documentstyle[11pt,eepic,leqno]{article}

\newtheorem{theorem}{Theorem}

\newtheorem{corollary}{Corollary}[theorem]

\newtheorem{proposition}{Proposition}[section]
\newtheorem{thm}[proposition]{Theorem}

\newtheorem{lemma}[proposition]{Lemma}

\newenvironment{remark}{\refstepcounter{proposition}\smallskip
\noindent {\bf Remark~\theproposition}\ }{\hskip\hsize plus0pt minus\hsize
\hbox{$\Box$}\smallskip}

\newenvironment{definition}{\refstepcounter{proposition}\smallskip
   \noindent{\bf Definition \theproposition}\ }{\smallskip}

\newenvironment{problem}{\refstepcounter{proposition}\smallskip
   \noindent{\bf Question \theproposition}\ }{\smallskip}

\newenvironment{proof}{\noindent{\bf Proof} \ }{\QED}\smallskip
\newcommand\QED{\hfill\raisebox{1pt}{\framebox[2.5mm][r]{
$\phantom x$}}\bigskip}

\newlength{\uoft}
\newlength{\uofe}

\title{A Momentum Construction for \\
Circle-Invariant K\"ahler Metrics\thanks{AMS Subject Classification:
32L05, 53C21 (Primary); 53C55, 32L07 (Secondary).}}
\author{\parbox{\uoft}{\center Andrew D. Hwang\thanks{Supported 
         in part by JSPS Fellowship \#P-94016 and an NSERC Canada Research
         Grant} \\
         University of Toronto}
\makebox[.6in]{}
 \parbox{\uofe}{\center Michael A. Singer \\
         University of Edinburgh}}
\date{}

\def \a#1 {\alpha_{[#1]}}
\def \A#1 {A_{[#1]}}
\def \Beta {{\rm B}}
\def \C {{\bf C}}
\def \Ex {E^{\times}}
\def \e {{\bf e}}

\def \lx {L^\times}

\def \O {{\cal O}}
\def \P {{\bf P}}
\def \R {{\bf R}}
\def \w {{\rm w}}
\def \Z {{\bf Z}}

\def \lapl {\lower 1pt\hbox{$\Box$}}
\def \aut {{\rm Aut}\,}
\def \dvol {{\rm dvol}_g}
\def \eps {\varepsilon}
\def \del {\partial}
\def \delbar {\bar\partial}
\def \ddbar {{\sqrt{-1}}\del\delbar}
\def \Im {{\rm Im}\,}
\def \Re {{\rm Re}\,}
\def \mom {{\mu}} %% Moment map as function of t
\def \mum {\tau}  %% Moment map; parameter on I
\def \Neg {\phantom{-}}
\def \tr {{\rm tr}}

\begin{document}

\maketitle
\tableofcontents

\section{Introduction}

The subject of this paper is the explicit construction of complete
K\"ahler metrics with prescribed---usually constant---scalar
curvature.  The technique, hereafter referred to as the {\em momentum
construction}, is a combination of two main ideas.  The first, which
goes back at least to the work of Calabi~\cite{Cal0}, is that of
constructing K\"ahler forms from K\"ahler potentials that are
essentially functions of one real variable.  The prototypical example
is the ansatz $\omega = \ddbar f(t)$, where
$t=(1/2)\log(|z_1|^2+\cdots+|z_n|^2)$, the $z_j$ are standard complex
coordinates on $\C^n$, and $f$ is a smooth function of one real
variable.  As is well-known, if $f(t)=t$, then $\omega$~descends to
the Fubini-Study metric on $\P^{n-1}$, while if $f(t) = e^t$, $\omega$
extends to the standard flat K\"ahler form on $\C^n$.

The general setting for a construction of this type is the total space
of an Hermitian holomorphic line bundle $p:(L,h)\to(M,\omega_M)$ over
a K\"ahler manifold. Let $t$ be the logarithm of the fibrewise
norm function defined by $h$ and consider the {\em Calabi
ansatz}
\begin{equation} \label{calans}
\omega = p^*\omega_M + 2 \ddbar f(t),
\end{equation}
which associates a K\"ahler form to (suitably convex) functions of one
real variable.

The second idea comes from symplectic geometry and involves a change
of variables $f(t)$ to $\varphi(\mum)$ closely related to the Legendre
transform.  For each $f$, (\ref{calans}) is invariant under the
circle-action that rotates the fibres of $L$. Let $X$ be the generator
of this action, normalized so that $\exp(2\pi X) =1$. Denote by $\tau$
the corresponding {\em moment map } determined (up to an additive
constant) by
\begin{equation} \label{intro2}
i_X\omega = - d\tau.
\end{equation}
(In fact $\tau = f'(t)$, see Section~2.1.)
%% The image of $\tau$ is some interval $I \subset \R$; we shall also
%% think of $\mum$ as an affine parameter on this interval. 
When $\omega$ comes from the Calabi ansatz, the function
$\|X\|^2_\omega$ is constant on the level sets of $\tau$, so there is
a function $\varphi:I\to(0,\infty)$---to be called the {\em momentum
profile} of~$\omega$---such that
\begin{equation} \label{intro3}
\varphi(\tau) = \|X\|^2_\omega.
\end{equation}
The crucial point is that $\omega$ can be reconstructed explicitly
from its profile; indeed $t$ and $\tau$ are related by the Legendre
transform, and the Legendre dual $F$ of $f$ satisfies
$F'' = 1/\varphi$.

The description of $\omega$ in terms of $\varphi$ has many advantages.
At one level, this is to be expected: (\ref{intro3}) shows that
$\varphi$ is a canonical geometric quantity, while $f$ determines
$\omega$ only through its second derivative. In particular, the only
condition that is needed for $\varphi$ to determine a K\"ahler metric
is that it be positive on (the interior of) $I$, cf.\
(\ref{intro3}). Positivity of (\ref{calans}), by contrast, corresponds
to two conditions on the derivatives of $f$. Further, the geometry of
the metric---for example its completeness or extendability properties
near the fixed-point set of the $S^1$-action---are easily read off
from the behaviour of $\varphi$ near the endpoints of~$I$.  However,
the decisive and most remarkable advantage, as far as the problem of
prescribed scalar curvature is concerned, is that the scalar curvature
of $\omega$ is given by a {\em linear} second-order differential
expression in $\varphi(\mum)$, in contrast to the fully nonlinear
fourth-order function that arises in the $f(t)$-description.
Consequently, the profiles that give rise to metrics of constant scalar
curvature are explicit rational functions of~$\mum$.

This dramatic simplification is already apparent in the basic example
of $S^1$-invariant metrics on $S^1$-invariant domains of $\P^1$. In
this case, the scalar curvature is given by $-\varphi''(\mum)/2$. For
the isometric $S^1$-action on the unit sphere in $\R^3$ by rotations
about the $z$-axis, the coordinate~$z$ itself is the moment
map~$\mum$, and $\varphi(\mum) = 1 -\mum^2$. It is immediately {\em
calculated}\/ that this metric has constant scalar curvature.
Generally, it is a pleasant and instructive exercise to classify
circle-invariant metrics with constant scalar curvature on subsets of
$S^2$ from this point of view, and to compare the simplicity of the
results with other approaches to this problem.  A sketch is provided
in Section~\ref{section:complete} (Geometry of fibre metrics; see also
Table~\ref{table-big}, page~\pageref{page-big}).
%% which contains a table of the results of this exercise.

This paper is organized as follows.  The remainder of the introduction
summarizes our main results. Section~\ref{section:KS} gives a
self-contained account of the momentum construction and
Section~\ref{section:inf-vol} applies the method to give general
existence theorems for complete K\"ahler metrics of constant scalar
curvature.  Finally Section~\ref{section:examples} is devoted to a
discussion of three topics: the scope and limitations of the momentum
construction; examples of line bundles $(L,h)$ 
%% on which complete K\"ahler metrics of constant scalar curvature can
%% be found---i.e.\ examples 
to which the results of Section~\ref{section:inf-vol} apply; and an
account of the related literature.  While postponing until that
section a careful explanation of the ways in which our work builds on
and extends that of previous authors, let us pause here to acknowledge
the most important sources of inspiration for the present work: the
papers of Calabi~\cite{Cal0}, Koiso--Sakane~\cite{KS1},
LeBrun~\cite{Lebrun}, and Pedersen--Poon~\cite{PP}.
\medskip

\noindent {\bf Acknowledgements}\quad The first-named author was
supported in part by a JSPS Postdoctoral Research Fellowship during
the early stages of this work, and by an NSERC Canada Individual
Research Fellowship. He is especially grateful to Professors Y.~Sakane
and T.~Mabuchi, and to Osaka University, for kind and extensive
hospitality, and to Professor J.~Bland, of the University of Toronto,
for many patient and helpful discussions.  The second-named author is
an EPSRC advanced fellow. He is also grateful to the University of
Toronto for its hospitality while some of this work was being carried
out.

\subsubsection*{Results about the momentum construction}

In the rest of the paper, $p:(L,h) \to (M,\omega_M)$ is an Hermitian
holomorphic line bundle, with curvature form $\gamma=-\ddbar\log h$,
over a K\"ahler manifold of complex dimension~$m$. A {\em compatible
momentum interval} is an interval $I \subset\R$, such that the closed
$(1,1)$-form $\omega_M(\mum):=\omega_M-\mum\gamma$ is positive for
every $\mum\in I$. The associated K\"ahler metric is
denoted~$g_M(\mum)$.

\begin{definition} 
\label{def:horizontal}
{\em Horizontal data} $\{p:(L,h)\to(M,g_M),\ I\}$ consists of an
Hermitian holomorphic line bundle over a K\"ahler manifold, together
with a compatible momentum interval. The whole assemblage is often
denoted $\{p,I\}$ for brevity.  A {\em momentum profile} is a function
$\varphi$ that is smooth on the closure of $I$ and positive on the
interior of $I$.
\end{definition}

The {\em completion} of a line bundle~$L$ is the $\P^1$ bundle
$\widehat{L}=\P(\O\oplus L)$, containing~$L$ as a Zariski-open subset
and obtained by adding a copy of~$M$ `at~$\infty$.' Let
$r:\widehat{L}\to[0,\infty]$ denote the continuous extension of the
square of the Hermitian norm function. All the metrics constructed
in this paper live on subsets of~$\widehat{L}$ obtained by
restricting~$r$ to an interval (or on manifolds obtained by partially
collapsing the zero and/or infinity sections).

\begin{definition}
\label{def:invariant}
Let ${\rm J} \subset [0,\infty]$ be an open interval. The
corresponding {\em invariant subbundle} $L'\subset \widehat{L}$ is the
$S^1$-invariant domain $r^{-1}({\rm J})$.
\end{definition}

Different choices of $\rm J$ yield six distinct complex-analytic fibre
types: ${\rm J}=[0,\infty]$ (the projective line); ${\rm J}
=[0,\infty)$ (the complex line); ${\rm J}=(0,\infty)$ (the punctured
line); ${\rm J}= [0,1)$ (the disk); ${\rm J}=(0,1)$ (the punctured
disk); and ${\rm J} = (e^{-l},e^l)$ (annuli).  In the last three
cases, homotheties have been used to reduce $\rm J$ to a standard
form. The invariant subbundles corresponding to the first five cases
will be denoted $\widehat{L}$, $L$, $L^\times$, $\Delta(L)$ and
$\Delta^\times(L)$ respectively.  Annulus-subbundles will not play a
major role in what follows.

\begin{definition} 
Let $L'\subset \widehat{L}$ be an invariant subbundle. A {\em
bundle-adapted metric} on $L'$ is a K\"ahler metric~$g$ whose K\"ahler
form~$\omega$ arises from the Calabi ansatz (\ref{calans}).
\label{def:adapt}
\end{definition}

The heart of the momentum construction is the fact, implicitly due to
Calabi and Koiso-Sakane, that if horizontal
data are given, then each momentum profile determines a unique
isometry class of bundle-adapted K\"ahler metric enjoying the
geometric properties of equations (\ref{calans})--(\ref{intro3}):

\begin{proposition} \label{prop:cons}
Let horizontal data $\{p,I\}$ and a momentum profile $\varphi$ be
given.  Then there exists an invariant subbundle $L'\subset
\widehat{L}$, unique up to homothety, and a bundle-adapted K\"ahler
metric~$g_\varphi$ on $L'$, unique up to isometry, with the following
properties: 
\begin{description}
\item{\rm(i)} The K\"ahler form $\omega_\varphi$ of $g_\varphi$ arises
from the Calabi ansatz; 
\item{\rm(ii)} The image of the moment map $\tau$ is~$I$;
\item{\rm(iii)} The normalized generator $X$ of the circle action
satisfies $g_\varphi(X,X)= \varphi(\tau)$. 
\end{description}
\end{proposition}

It is natural to ask how the geometric properties, especially
completeness and the scalar curvature~$\sigma_\varphi$, of~$g_\varphi$
are encoded by~$\varphi$. As will be seen in 
Section~\ref{section:complete} below, completeness of~$g_\varphi$ is
encoded by the boundary behaviour of~$\varphi$, in fact by the 2-jet
at the endpoints of~$I$. As claimed above, the scalar curvature is
linear in~$\varphi$:

\begin{theorem}
\label{thm:CKS}
Let horizontal data be given. For each $\mum\in I$, let
$\sigma_M(\mum)$ denote the scalar curvature of $g_M(\mum)$, and
define $Q:I\times M\to\R$ by 
$Q(\tau) = \omega_M(\tau)^m/\omega_M^m$.  Then the scalar curvature of
$g_\varphi$ is given by 
\begin{equation}
\label{eqn:CKS2}
\sigma_\varphi = \sigma_M(\mum)-{1\over{2Q}} 
{{\del^2}\over{\del\mum^2}}(Q\varphi)(\mum).
\end{equation}
\end{theorem}
To interpret~(\ref{eqn:CKS2}), regard the $S^1$-invariant function
$\sigma_\varphi$ on~$L'$ as a function on $I\times M$ by factoring
through the $S^1$-action.

The improvement of Theorem~\ref{thm:CKS} over earlier results is that
no curvature hypotheses are imposed on the horizontal data.  Of
course, it is too much to expect that for arbitrary horizontal data,
$\sigma_\varphi$ can be made constant by an appropriate choice of
profile (a single function of one variable).  Natural {\em
sufficient}\/ curvature hypotheses are to assume the two terms
in~(\ref{eqn:CKS2}) separately depend only on~$\mum$ for every
profile: 

\begin{definition}  Horizontal data are said to be {\em
$\sigma$-constant}\/ if
\begin{description}
\item{\rm(i)} The curvature endomorphism $\Beta=\omega_M^{-1}\gamma$
has constant eigenvalues on $M$; 
\item{\rm(ii)} The metric $g_M(\mum)$ has constant scalar curvature
for each~$\mum\in I$. 
\end{description} \label{sconst}
\end{definition}

The simplest examples of $\sigma$-constant horizontal data are
pluricanonical bundles over Einstein-K\"ahler manifolds. More
generally, if $g_M$ has constant scalar curvature, and $\gamma(L,h)$
is a multiple of $\omega_M$, then the data $p:(L,h) \to (M,\omega_M)$
are $\sigma$-constant. More specific examples appear in
Section~\ref{section:examples}. 

The conclusion of Theorem~\ref{thm:CKS} is due to  Guan~\cite{Guan1}
and Hwang~\cite{Hwang} under stronger curvature hypotheses
(``$\rho$-constancy,'' see Section~\ref{section:compare}).  The
corresponding statement for $\sigma$-constant data is used repeatedly
in the sequel:
\begin{corollary}
\label{cor:CKS}
Let the horizontal data $\{p:L\to M,\ I\}$ be $\sigma$-constant. Then
$\sigma_M(\tau) = P/2Q$, where $Q$ is as in Theorem~\ref{thm:CKS},
$P$ is a polynomial in $\tau$, and for each momentum profile~$\varphi$,
\begin{equation}
\label{eqn:CKS3}
\sigma_\varphi={1\over{2Q}}\Bigl(P-(Q\varphi)''\Bigr)(\mum).
\end{equation}
\end{corollary}

\subsubsection*{Metrics of constant scalar curvature on line bundles}

The next two theorems are generalizations of the work
of various authors. Detailed attribution is given in
Section~\ref{section:compare}.  

\begin{theorem} \label{thm:A} 
Let $I=[0,\infty)$, and let $\{p,I\}$ be $\sigma$-constant horizontal
data with $\gamma \leq 0$, $\gamma\not= 0$, and with $g_M$~complete.
Then there exists a real number $c_0$ with the following property: For
every $c < c_0$, and for at most finitely many $c>c_0$, the disk
bundle~$\Delta(L)$ carries a complete K\"ahler metric~$g_c$, of scalar
curvature~$c$, whose restriction to the zero section is~$g_M$.  When
$c=c_0$, the analogous conclusion holds, but the metric lives on the
total space of~$L$.  For $c\leq0$, the metric~$g_c$ is Einstein iff
${{c}\over{m+1}}\omega_M=\rho_M+\gamma$.
\end{theorem}

An analogue of Theorem~\ref{thm:A} holds for punctured disk
bundles. Of course, the metrics are not obtained from the metrics in
Theorem~\ref{thm:A} by restriction, since removing the zero section
leaves an incomplete metric.  

\begin{theorem} \label{thm:B} 
Let $I=(0,\infty)$, and let $\{p,I\}$ be $\sigma$-constant horizontal
data with $\gamma \leq 0$, $\gamma\not= 0$, and with $g_M$~complete.
Then there exists a real number $c_0^\times$ with the following
property: For every $c < c_0^\times$, and for at most finitely many
$c>c_0^\times$, $\Delta^\times(L)$ carries a complete K\"ahler
metric~$g_c^\times$, of scalar curvature~$c$, whose symplectic
reduction at $\mum=1$ is~$g_M-\gamma$.  When $c=c_0^\times$, the
analogous conclusion holds, but the metric lives on the total space
of~$\lx$.  For $c<0$, the metric~$g_c$ is Einstein iff
${{c}\over{m+1}}\omega_M=\rho_M$.
\end{theorem}

\begin{remark} There is some redundancy in these statements;
the condition $I=(0,\infty)$ implies that $\gamma$ is non-positive,
for $\omega_M(\tau)$ is supposed to be positive for all $\tau \in
I$.  On the other hand, it is natural to exclude the case $\gamma=0$,
since bundle-adapted metrics on flat bundles are local product
metrics, and can thus be understood in an elementary fashion.
\end{remark}

\begin{remark} 
The choice $I=[0,\infty)$ or $(0,\infty)$ is a normalization of the
data which results in no loss of generality, see
Lemma~\ref{lemma:interval} below.
\end{remark}

\begin{remark}
The constants $c_0$ and $c_0^\times$ are roots of polynomials in
one variable and can be estimated in terms of the horizontal data. In
good cases they can  be found exactly.  In particular, if $\gamma$ is
negative-definite then $c_0=0$ and one has the pleasant conclusion
that for every $c<0$, $\Delta(L)$ admits a complete K\"ahler metric
with scalar curvature~$c$, while $L$ admits a complete scalar-flat
K\"ahler metric.
\end{remark}

Theorem~\ref{thm:A} contains many previously-known results (see
Section~\ref{section:examples}), but even when $M$ is a complex
curve, an interesting new result is obtained:
\addtocounter{theorem}{-1}
\begin{corollary}
\label{cor:flat-surf} 
Let $M= \C$, $\omega_M =\frac{\sqrt{-1}}{2}dz\wedge d\bar{z}$ the
standard flat K\"ahler form, $(L,h)\to M$ the trivial line bundle
equipped with an Hermitian metric of curvature $\gamma=-2\pi
k\,\omega_M$, {\rm(}$k$ a positive constant{\rm)}. Then the total
space of $L$---biholomorphic to $\C^2$---admits a complete,
scalar-flat K\"ahler metric that is not Ricci-flat, and whose fibre
metric is asymptotically cylindrical.
\end{corollary}

The importance of this example is that scalar-flat K\"ahler metrics on
complex surfaces are anti-self-dual (in the sense of $4$-dimensional
conformal geometry), and are thus of independent interest.

The following Einstein-K\"ahler metrics arising from
Theorem~\ref{thm:A} are due to Calabi~\cite{Cal0} when the base is
Einstein-K\"ahler (i.e.\ when $n=1$ in Corollary~\ref{cor:A1} below).
Moreover, non-homothetic metrics are obtained in
Corollary~\ref{cor:A1} by scaling the $\omega_{M_i}$ separately.  In
Corollary~\ref{cor:A2}, further metrics arise by tensoring with a flat
line bundle when the base is not simply-connected.

\begin{corollary}
\label{cor:A1}
For $1\leq i\leq n$, let $p_i:K_i\to M_i$ be the canonical bundle of a
compact Einstein-K\"ahler manifold of positive curvature. Then there
is a complete, Ricci-flat K\"ahler metric on the total space
of~$\bigotimes_i p_i^*K_i$.
\end{corollary}

\begin{corollary}
\label{cor:A2}
For $1\leq i\leq n$, let $p_i:K_i\to M_i$ be a compact
Einstein-K\"ahler manifold of positive curvature. Then for all
$\ell>1$, there is a complete Einstein-K\"ahler metric of negative
curvature on the disk subbundle of~$\left(\bigotimes_i
p_i^*K_i\right)^\ell$. If $M_i$ are complete, negative
Einstein-K\"ahler manifolds, then the disk subbundle in $(\bigotimes_i
p_i^* K_i^{-1})^\ell$ admits a complete, negative Einstein-K\"ahler
metric.
\end{corollary}
\addtocounter{theorem}{1}
\setcounter{corollary}{0}

Theorem~\ref{thm:B} yields some new Einstein-K\"ahler metrics on
punctured disk bundles. The fibre metric has one end of finite area
and one end of infinite area.

\begin{corollary}
\label{cor:B1} 
Let $p:L\to(M,g_M)$ be a holomorphic line bundle over a complete,
negative Einstein-K\"ahler manifold, and assume $c_1(L)\in
H^{1,1}(M,\Z)$ is represented by a form~$\gamma$ whose eigenvalues
with respect to~$\omega_M$ are constant on~$M$. Then $\Delta^\times(L)$
admits a complete, negative Einstein-K\"ahler metric.
\end{corollary}

\subsubsection*{Metrics of constant scalar curvature on vector bundles}

The Calabi ansatz applies, in modified form, to vector bundles~$E\to
D$ of rank $n>1$; Calabi's construction of a complete, Ricci-flat
metric on the cotangent bundle of a compact, rank-one Hermitian
symmetric space~\cite{Cal0} is the prototypical example. 
%Unfortunately, the ansatz does not usually
%reduce the problem to an ODE, even if $D$ is Einstein--K\"ahler and
%$E$ has a Hermitian-Einstein metric; Calabi's construction relies on
%homogeneity of the base and irreducibility of~$E$. 
In seeking generalizations of this example, one adopts the following
strategy. Take $M = \P(E)$, $L =\tau_E$ the tautological line bundle
over $\P(E)$.  Since the total space of $L$ is the blow-up along the
zero-section of $E$, a metric on $E$ is the same thing as a suitably
degenerate metric on $\tau_E$.  The momentum construction is easily
modified to produce such partially degenerate metrics on~$L$: 
\begin{theorem}
\label{thm:C}
Let $(E,h)\to(D,g_D)$ be an Hermitian holomorphic vector bundle of
rank~$n$ over a complete K\"ahler manifold of dimension~$d$, and let
$p:(\tau_E,h)\to M=\P(E)$ be the tautological bundle with the
induced Hermitian structure. Let $\gamma$ denote the curvature form of
$(\tau_E,h)$.  Assume $\omega_M(\mum):=\omega_D-\mum\gamma$ is
a K\"ahler form on~$M$ for all $\mum>0$, and that
$p:(\tau_E,h)\to\bigl(M,g_M(\mum_0)\bigr)$ is $\sigma$-constant for
some~$\mum_0>0$.  Then there is a $c_0$ such that for every $c <
c_0$, and for at most finitely many $c>c_0$, the momentum construction
determines a complete K\"ahler metric with scalar curvature~$c$ on the
disk subbundle~$\Delta\subset E$.  When $c\leq0$, this metric is
Einstein if and only if ${c\over{n+d}}\omega_D=\rho_M+n\gamma$.
Furthermore, there is a complete K\"ahler metric on $E$ with
scalar curvature~$c_0$.
\end{theorem}

Unfortunately, $\sigma$-constancy of the data $\tau_E \to \P(E)$ is a
strong assumption.  It is not even sufficient for $E$ to be a
homogeneous vector bundle, as the example $E=\O\oplus\O(k)\to\P^1$
($k<0$) shows.  In this case $\P(E)$ is a Hirzebruch surface, which
does not admit a K\"ahler metric of constant scalar curvature, so
condition~(ii) in the definition of $\sigma$-constancy is never
satisfied.

Each of the next three corollaries contains tractable hypotheses that
guarantee the applicability of Theorem~\ref{thm:C}.  In
Corollary~\ref{cor:homog}, it is assumed that both $E$~and $\P(E)$ are
homogeneous. The pseudoconvexity hypothesis guarantees a metric of
non-positive curvature on the tautological bundle.  The hypotheses of
Corollary~\ref{cor:sum} imply that $\P(E)$ is a product and that
$\tau_E$ is a tensor product of line bundles pulled back from the
factors.  Corollary~\ref{cor:stable}, which rests on the result of
Narasimhan-Seshadri to the effect that a stable bundle over a curve is
projectively flat, gives families of examples depending on continuous
parameters in the event that the base of~$E$ is one-dimensional.
As above, the restriction of the advertized metric to the zero section
is the base metric~$g_D$. In particular, varying the scalar curvature
does not vary the metric through homotheties.

\begin{corollary}
\label{cor:homog}
Let $E\to D$ be a homogeneous vector bundle of rank~$n>1$ over a
compact, homogeneous K\"ahlerian manifold. Assume further that the
projective bundle~$M=\P(E)$ is homogeneous, and that the zero section
of~$E$ has a pseudoconvex tubular neighborhood.  Endow the base with a
homogeneous K\"ahler metric~$g_D$, and the bundle with a homogeneous
Hermitian metric~$h$. Then for every $c<0$, there is a complete
K\"ahler metric~$g_c$ of scalar curvature~$c$ on the disk
bundle~$\Delta$, and there is a complete, scalar-flat metric on the
total space of~$E$.
\end{corollary}

\begin{corollary}
\label{cor:sum}
Suppose $(D,g_D)$ is complete and that $(\Lambda,h)\to(D,g_D)$ is
$\sigma$-constant, with curvature $\gamma<0$. Then the total space of
$E=\Lambda\otimes\C^n=\Lambda^{\oplus n}$ admits a complete,
scalar-flat metric, and for every $c<0$, the disk bundle~$\Delta$
{\rm(}taken with respect to the natural Hermitian structure{\rm)}
admits a complete K\"ahler metric of scalar curvature~$c$.
\end{corollary}

\begin{corollary}
\label{cor:stable}
Let $p:E\to C$ be a stable vector bundle of rank~$n$ and
degree~$k\leq0$ over a compact Riemann surface~$C$, of genus $g\geq 2$
and endowed with the unit-area metric of constant Gaussian
curvature. Then the total space of~$E$ admits a complete, scalar-flat
K\"ahler metric, and there exists an Hermitian structure~$h$ such that
for every $c<0$, the disk subbundle of $(E,h)$ carries a complete
K\"ahler metric of scalar curvature~$c$.  Finally, if $k\leq2-2g$,
then there is a complete Einstein-K\"ahler metric with scalar
curvature $c=2\pi(n+1)(2g-2+k)$, on~$E$ if $c=0$ and on $\Delta\subset
E$ if $c<0$.
\end{corollary}

\subsubsection*{Limitations of the Calabi ansatz}

If horizontal data are $\sigma$-constant, then by~(i) of
Definition~\ref{sconst}, $Q(\mum)$ is constant for all $\mum\in I$, so
the second term in~(\ref{eqn:CKS2}) depends only on~$\mum$ for every
profile~$\varphi$, while (ii)~says exactly that $\sigma_M(\mum)$ is
constant for all $\mum$.  Conversely, it is not difficult to show that
if, {\em for every}\/ $\varphi$, $\sigma_\varphi$ depends only upon
$\mum$ then the horizontal data are $\sigma$-constant.  We believe the
same conclusion follows merely if {\em there exists}\/ a profile
inducing a metric of constant scalar curvature, but are at present
able to prove only a partial result in this direction:

\begin{theorem}
\label{thm:constant} 
Let $\{p,I\}$ be horizontal data with $I=[0,\eps)$ for some $\eps>0$,
and assume there exist real-analytic functions $\sigma_1$~and
$\sigma_2$ on~$I$, and distinct profiles $\varphi_1$~and $\varphi_2$,
with $\varphi_1(0)=\varphi_2(0)$ and inducing metrics of scalar
curvature $\sigma_1(\mum)$~and $\sigma_2(\mum)$, respectively. Then
the horizontal data are $\sigma$-constant. 
\end{theorem}

Theorem~\ref{thm:constant} is proved by considering the difference
$\varphi_1-\varphi_2$, and is of course purely local.  A simple but
suggestive corollary is that on horizontal data that are not
$\sigma$-constant, there is at most one metric of constant scalar
curvature arising from the Calabi ansatz. In the statement below, it
is not necessary to assume $c_1 \not= c_2$; if $c_1=c_2$, then
$\varphi_1$~and $\varphi_2$ are distinct iff
$\varphi_1'(0)\not=\varphi'_2(0)$.

\begin{corollary}
\label{cor:constant}
Let $\{p,I\}$ be horizontal data as above, and assume there exist
distinct profiles $\varphi_1$~and $\varphi_2$, equal at~$0$ and
inducing metrics of scalar curvature $c_1$~and $c_2$,
respectively. Then the horizontal data are $\sigma$-constant.
\end{corollary}

There does not seem to be an easy way to prove $\sigma$-constancy by
assuming existence of {\em one}\/ profile, without additional
hypotheses. Further partial results and a conjectural strengthening of
Theorem~\ref{thm:constant} are described in
Section~\ref{section:examples}.

\section{Bundle-Adapted Metrics}
\label{section:KS}
\setcounter{equation}{0}

This section gives a full and self-contained account of the momentum
construction. To facilitate the presentation, a detailed account of
the credit due to earlier authors is postponed to
Section~\ref{section:compare}.  We begin with a short summary of the
most important formulae of the momentum construction.  Proofs and
further comments come in Section~\ref{section:complete}.

\subsection{The Momentum Construction}

Terminology is as in Definitions
\ref{def:horizontal}--\ref{def:adapt}. The following notational 
conventions are used systematically: A K\"ahler metric is
denoted~$g$. Its K\"ahler form, Ricci form, scalar curvature, and
Laplace operator (acting on functions) are denoted $\omega$, $\rho$,
$\sigma$, and $\lapl$ respectively. Subscripts are used
descriptively when several metrics are under consideration.

Fix $b\leq\infty$. Let~$I\subset\R$ be an interval with interior
$(0,b)$, and pick $\mum_0 \in (0,b)$. For each profile~$\varphi$,
define 
\begin{equation} \label{t12def}
t_1=\lim_{\mum\to0^+}\int_{\mum_0}^\mum {{dx}\over{\varphi(x)}},\quad
t_2=\lim_{\mum\to b^-}\int_{\mum_0}^\mum {{dx}\over{\varphi(x)}},
\end{equation}
noting that $-\infty\leq t_1<t_2\leq\infty$. Introduce functions
$\mom$, $s$, and $f$ of $t$ by 
\begin{equation}
\label{eqn:mu}
t  = \int_{\mum_0}^{\mom(t)} {{dx}\over{\varphi(x)}},\qquad 
s(t) = \int_{\mum_0}^{\mom(t)} {{dx}\over{\sqrt{\varphi(x)}}},\qquad
f(t) = \int_{\mum_0}^{\mom(t)} {{x\,dx}\over{\varphi(x)}}. 
\end{equation}
Then the equations
\begin{equation}
\label{eqn:tau-ode}
\mom'=\varphi\circ\mom,\qquad
s'=\sqrt{\varphi\circ\mom},\qquad
f'=\mom
\end{equation}
hold on $(t_1,t_2)$.

Suppose now that $\{p:(L,h)\to(M,g_M),I\}$ are horizontal data, and
let $t: \lx \to \R$ be the logarithm of the norm-function. The
functions $\mom$, $s$, and $f$ of (\ref{eqn:mu}) induce respective
functions $\mum:=\mom(t)$, $s(t)$, and $f(t)$ on $L' = \{z\in
L^\times: t_1 < t(z) <t_2\}$, and determine a closed $(1,1)$-form
\begin{equation}
\label{eqn:omega}
\omega_\varphi=p^*\omega_M+dd^c f(t)=p^*\omega_M+2\ddbar f(t).
\end{equation}
It transpires that $\omega_\varphi$ is positive, i.e.\ is a K\"ahler
form. The associated bundle-adapted metric will be denoted
$g_\varphi$, and the mapping $\varphi \longmapsto \omega_\varphi$ will
be called the {\em momentum construction}.

\subsubsection*{The fibre metric}

If $z^0$ is a linear coordinate on a fibre of $L$ then the restriction
of $g_\varphi$ to the fibre is given by
\begin{equation} \label{fibg}
g_{\rm fibre} = 
\varphi(\mum)\,\left|\frac{dz^0}{z^0}\right|^2.
\end{equation}
Thus $\varphi(\mum)$ is interpreted as the conformal factor relating
$g_{\rm fibre}$ to the flat metric on the cylinder.  The function~$s$
of~(\ref{eqn:mu}) is the geodesic distance for $g_{\rm fibre}$, while
$2\pi\bigl(\mom(b)-\mom(a)\bigr)$ is the area of the subset of the
fibre given by $a \leq t \leq b$.

\subsubsection*{$S^1$-invariant functions on $L$}

Fix a bundle-adapted metric on~$L'$.  The functions $t$~and
$\mum=\mom(t)$ have the same level sets on~$L'$, but while $t$ depends
only on the Hermitian structure of~$L$, $\mum$ depends on the profile
as well.  By a customary abuse of notation, $t$~and $\mum$ are
regarded as variables in the intervals $(t_1,t_2)$ and $I$,
respectively.

Define $\pi: L' \to I\times M$ by $\pi=(\mum,p)$. Each
fibre of $\pi$ is an orbit of the $S^1$-action, so circle-invariant
tensors on $L'$ may be identified with tensors on $I\times M$; this
will henceforth be done freely, with pullback~$\pi^*$ suppressed.
These identifications are diagrammed in Figure~\ref{fig:cscmaps}~(a);
the vertical maps on the right are projections. The functions $\mom$,
$s$, and $f$ of equation~(\ref{eqn:mu}) are defined on $(t_1,t_2)$,
while $\varphi$ and other functions of geometric interest are defined
on~$I$.

\begin{figure}[hbt] 
\label{fig:cscmaps}
%% Figures 2.1 (a) and (b); each sits in a minipage, and consists of a
%% picture and centered caption. The entire figure consists of one
%% centered minipage, containing the two sub-minipages.

\begin{minipage}[b]{\textwidth}
\begin{center}

	\begin{minipage}[b]{192pt}
	\setlength{\unitlength}{12pt} %% \qquad

	\begin{picture}(8,7)(-2,0)

		%% Spaces
		\put(0,.5){$(t_1,t_2)$}
		\put(.5,4){$L'$}
		\put(10.75,.5){$I$}
		\put(10,4){$I\times M$}
		\put(10.5,7.5){$M$}

		%% Maps
		\put(1.25,2){$t$}
		\put(6,1){$\mu$}
		\put(6,2.5){$\mum$}
		\put(6,4.5){$\pi$}
		\put(6,6.5){$p$}

		%% Arrows
		\put(1,3.5){\vector(0,-1){2}}
		\put(11,3.5){\vector(0,-1){2}}
		\put(11,5){\vector(0,1){2}}
		\put(2,4.25){\vector(1,0){7.75}}
		\put(2,4.75){\vector(3,1){8}}
		\put(2,3.75){\vector(3,-1){8}}
		\put(2.5,.75){\vector(1,0){7.5}}
	\end{picture}

	\begin{center}
	Figure~2.1 (a): Relation between maps in the
	momentum construction
	\end{center}
	\end{minipage}
\hfil
	\begin{minipage}[b]{3in}
	\setlength{\unitlength}{.25in}
		
	\begin{picture}(12,6.5)(0,.75)
		{\thinlines
		\put(1.667,7){\line(1,0){5.667}}
		\put(1.667,7){\line(1,-3){1.333}}
		\put(6,3){\line(1,3){1.333}}
		\put(4.5,2.5){\vector(0,-1){1}}
		\put(7.5,4.5){\vector(1,0){3}}
		\put(4,3){\line(-1,6){.667}}
		\put(5,3){\line(1,6){.667}} }

		{\thicklines
		\put(3,1){\line(1,0){3}}
		\put(2,6){\line(1,0){5}}
		\put(11,3){\line(0,1){4}}
		\put(3,3){\line(1,0){3}} }

		\put(.5,4){$I\times M$}
		\put(7.75,6){$\omega_M(\mum)$}
		\put(6.5,3){$\omega_M$}
		\put(4.25,.25){$M$}
		\put(12,4.5){$I$}
		\put(4.75,2){$p$}
		\put(8.75,4){$\mum$}

		\put(10.85,5.9){$\bullet$}
		\put(11.5,6){$\mum$}
		\put(10.85,2.9){$\bullet$}
		\put(11.5,3){$0$}
	\end{picture}

	\begin{center}
	Figure~2.1 (b): Symplectic reduction of~$g_\varphi$
	\end{center}

	\end{minipage} 
	\end{center}
\end{minipage}
\end{figure}

Two instances of identification by~$\pi$ deserve immediate mention.
First, the Euler vector field $-JX$ on $L'$ pushes forward to
$\varphi(\mum){\del\over{\del\mum}}$ on $I\times M$, cf.\
equation~(\ref{ddbtau}) below. The second example arises from the
1-parameter family $\omega_M-\mum\,\gamma$ of K\"ahler forms on~$M$,
regarded as a $(1,1)$-form $\omega_M(\mum)$ on $I\times M$, see
Figure~\ref{fig:cscmaps}~(b). The family of Ricci forms
$\rho_M(\mum)$, and scalar curvature functions $\sigma_M(\mum)$, are
similarly regarded as living on $I\times M$. {\em These tensors depend
only on the horizontal data}. If a profile is specified in addition,
then there is a map $\pi:L'\to I\times M$ as above, and each of these
tensors is identified with an $S^1$-invariant tensor on~$L'$.

\subsubsection*{The Laplacian, Ricci and scalar curvature}

\begin{definition} \label{def:pqr} 
The {\em curvature endomorphism}\/ $\Beta$ is $\omega_M^{-1}\gamma$,
the symmetric endomorphism of $T^{1,0}M$ obtained by raising the
second index of~$\gamma$. Similarly the {\em Ricci endomorphism}
$\varrho$ of $T^{1,0}M$ is defined by raising the second index of the
Ricci form~$\rho_M$.  Put
\begin{eqnarray}
Q(\mum) & = & {{\omega_M(\mum)^m}\over{\omega_M^m}} 
\quad=\quad\det({\rm I}-\mum \Beta) \label{Qdef} \\
R(\mum) & = & \tr_{\omega_M(\mum)}\rho_M
\quad=\quad\tr\Bigl[({\rm I}-\mum \Beta)^{-1}\varrho\Bigr]
\label{Rdef} \\
P(\mum) & = & 2QR(\mum)\vphantom{\Bigg|}\nonumber%% \label{Pdef}
\end{eqnarray}
%=\quad\sum_{\nu=1}^m
%{{\varrho_\nu}\over{1-\beta_\nu\mum}}\nonumber
%\end{eqnarray}
\end{definition}

\noindent Note that $P$ and $Q$ are smooth functions on $I\times M$, with
polynomial dependence upon $\mum$, while $R$ may be thought of as a
rational function on~$I$, with coefficients depending smoothly on~$z$
in~$M$. The notation reflects the fact that when the horizontal data
are $\sigma$-constant, these functions depend only on~$\mum$, i.e.\
are constant on~$M$ for every $\mum\in I$.

On a K\"ahler manifold, if $g_{i\bar{\jmath}}$ are the components of
the metric in local holomorphic coordinates, then the Ricci form is
given locally by $\rho = -\ddbar \log V$, where
$V=\det(g_{i\bar{\jmath}})$. For a bundle-adapted metric, there is a
choice of coordinates such that the quantities $V_\varphi$ and $V_M$
satisfy
\begin{equation} \label{volumef}
V_\varphi = (\varphi Q)(\mum)\, V_M.
\end{equation}
The Ricci form of $\omega_\varphi$ is therefore given by
\begin{equation} \label{ricc}
\rho_\varphi = p^*\rho_M - \ddbar\log \varphi Q(\mum),
\end{equation}
and the scalar curvature is found by taking the trace:
\begin{equation} \label{sigma1}
\sigma_\varphi = R(\mum) - \lapl_\varphi \log \varphi Q(\mum),
\end{equation}
$\lapl_\varphi$ being the $\delbar$-Laplacian of $\omega_\varphi$.
This Laplacian has a reasonably pleasant expression in terms of
$\lapl_{\omega_M(\mum)}$; for each smooth function~$\psi$ on $I\times
M$, 
\begin{equation} \label{boxphi}
\lapl_\varphi \psi = \lapl_{\omega_M(\mum)} \psi(\mum,\cdot)
+ \frac{1}{2Q}{\del\over{\del\mum}}\left[
\varphi Q(\mum)\, {{\del \psi}\over{\del\mum}}\right].
\end{equation}
Applying this to the function $\psi=\log(\varphi Q)$ and
combining with~(\ref{sigma1}) yields
\begin{equation} \label{sigma2}
\sigma_\varphi 
= R(\mum) - \lapl_{\omega_M(\mum)}\log  Q(\mum) 
- \frac{1}{2Q} \frac{\del^2}{\del \mum^2}\bigl(\varphi Q\bigr)(\mum),
\end{equation}
since $\lapl_{\omega_M(\mum)}\log \varphi Q(\mum) =
\lapl_{\omega_M(\mum)}\log  Q(\mum)$, $\varphi$ being independent of
$z\in M$.   
The first two terms together make up the scalar curvature
$\sigma_M(\mum)$ of $\omega_M(\mum)$, so 
\begin{equation} \label{sigma3}
\sigma_\varphi = \sigma_M(\mum) - \frac{1}{2Q}
\frac{\del^2}{\del \mum^2}\bigl(\varphi Q\bigr)(\mum).
\end{equation}
If the data are $\sigma$-constant, (\ref{ricc}) and
(\ref{sigma2}) simplify as follows:
\begin{eqnarray}
\rho_\varphi & = & p^*\rho_M + 
\frac{1}{2Q}(\varphi Q)'(\mum)\,p^*\gamma - \frac{1}{2\varphi}\left[
\frac{1}{Q}(\varphi Q)'\right]'(\mum)\, d\mum\wedge d^c\mum 
 \label{cricc} \\
\sigma_\varphi & = & R(\mum) -  \frac{1}{2Q}
\frac{\del^2}{\del \mum^2}\bigl(\varphi Q\bigr)(\mum). \label{scs}
\end{eqnarray}
These, together with the expression
\begin{equation} \label{ommu}
\omega_\varphi = p^*\omega_M - \mum\, p^*\gamma+{1\over{\varphi}}\,
d\mum\wedge d^c\mum,
\end{equation}
are the main formulae that will be needed in Section~3 in the
construction of metrics with prescribed scalar curvature by ODE
methods.

\subsection{Proofs}

We shall now work through the proofs of the statements in the summary
just given.  We begin by fixing some notation and recalling some basic
formulae.

\subsubsection*{Adapted coordinates}

Let $V\subset M$ be a coordinate chart over which~$L$ is trivial. Then
there exists a {\em line bundle chart}, namely, a local coordinate
system $z^0,z^1,\ldots,z^m$ for~$L$ in which $z^0=\rho e^{i\theta}$~is
a fibre coordinate and $z=(z^1,\ldots,z^m)$ are (pullbacks of)
coordinates on~$M$.  Greek indices run from $1$~to $m$, while Latin
indices run from~$0$ to~$m$.

In such a chart, there is a smooth, positive function~$h:V\to\R$ such
that $r=z^0\bar{z}^0\,h(z)$; under change of chart, the local
function~$h$ is multiplied by the norm squared of a non-vanishing
local holomorphic function in~$V$.  In a line bundle chart,
the Euler vector field is given by
$$\Upsilon=z^0{{\del}\over{\del z^0}}={1\over2}\left(
\rho{\del\over{\del\rho}}-\sqrt{-1}{\del\over{\del\theta}}\right),$$
while twice the imaginary and real parts are described variously as
$$X = -2\,\Im\Upsilon=\sqrt{-1}(\Upsilon-\overline\Upsilon)
={\del\over{\del\theta}},\qquad H = -JX =
2\,\Re\Upsilon=(\Upsilon+\overline\Upsilon)
=\rho{\del\over{\del\rho}}.$$ In particular this gives a local formula
for the normalized generator $X$ of the $S^1$ action on $L$. 
It is sometimes convenient to use the fibre
coordinate $w^0=\log z^0=:\zeta+i\theta$, in which case
$\Upsilon=\del/\del w^0$.

The level sets of~$r$ are real hypersurfaces in~$L$, and their tangent
spaces are the horizontal spaces of the Hermitian connection
$\Theta=\del\log r=2\del t$ of~$(L,h)$. For each point~$x$ of~$M$,
there exists a line bundle chart $(z^0,z)$ such that $z^\alpha(x)=0$
for $1\leq\alpha\leq m$ and $\del_\alpha r=0$ on the fibre~$L_x$. Such
a coordinate system is said to be {\em adapted} to~$(L,h)$ at~$x\in
M$. In adapted coordinates at~$x$, the connection form $\Theta$ is
equal to $dw^0=dz^0/z^0$ along the fibre~$L_x$.

\subsubsection*{Sign conventions and Levi forms}

The wedge product is normalized so that interior multiplication is a
(graded) derivation: If $V$ is a vector and $\xi$ and $\eta$ are
$1$-forms, then $i_V(\xi \wedge\eta) = \langle V,\xi \rangle \eta -
\langle V,\eta \rangle \xi$.  Extensive use will be made of the real
operators $d=\del+\delbar$ and $d^c=\sqrt{-1}(\delbar-\del)$, and of
the useful formulae
$$dd^c\,u = 2\ddbar u,\qquad 
du\wedge d^c u = 2\sqrt{-1}\del u \wedge \delbar u,$$ 
which hold for every smooth function~$u$. Finally if~$g$ has
components $(g_{i\bar{\jmath}})$ in local holomorphic
coordinates~$(z^i)$, then $\omega = 
\frac{\sqrt{-1}}{2}g_{i\bar{\jmath}}dz^i\wedge d\bar{z}^j$.

The curvature form~$\gamma$ of an Hermitian line bundle satisfies
$-dd^c t = - \ddbar \log r = p^*\gamma$. Combining this with the chain
rule gives the simple but important formula 
\begin{equation} \label{levi-t}
dd^c u(t) = u''(t)\,dt\wedge d^ct - u'(t)\,p^*\gamma
\end{equation}
when $u$ is a smooth function of one variable. This is often used in
the form
\begin{equation} \label{ddbtau}
dd^c\psi(\mum) = -\bigl(\varphi \psi'\bigr)(\mum)\,p^*\gamma 
+ {1\over\varphi}\bigl(\varphi \psi'\bigr)'(\mum)\,
d\mum\wedge d^c\mum,
\end{equation}
which follows from (\ref{levi-t}) and the first of (\ref{eqn:tau-ode}).

\subsubsection*{Proof of Proposition~\protect{\ref{prop:cons}}}

Referring to the momentum construction, especially formulae
(\ref{t12def})--(\ref{eqn:omega}), it must be shown that
$g_\varphi$ is a K\"ahler metric (i.e.\ that $\omega_\varphi>0$), and
that the moment map is given by~$\mum$. Expanding $dd^c f$ using
(\ref{ddbtau}) gives (\ref{ommu}): 
$$
\omega_\varphi 
= p^*\omega_M - \mum\, p^*\gamma + \varphi\, dt\wedge d^ct 
= p^*\omega_M - \mum\, p^*\gamma + 
{1\over\varphi}\, d\mum\wedge d^c\mum.
$$
Since $\omega_M -\mum\gamma$ and~$\varphi$ are both positive on~$I$ by
hypothesis, it follows that $\omega_\varphi$ is positive. Now
$$
i_X(dt) = 0 = i_X\Bigl(\omega_M(\mum)\Bigr),
$$
so $i_X\omega_\varphi = -\varphi\,dt \wedge i_X(d^ct) = -\varphi\,dt$
from the formula
$$
d^ct = \frac{\sqrt{-1}}{2}\left(\frac{d\bar{z}^0}{\bar{z}^0} -
\frac{d z^0}{z^0}\right) + \frac{1}{2}d^c\log h
$$ 
in a line-bundle chart. But $\varphi\,dt = d\mum$ by
(\ref{eqn:tau-ode}), so $\mum$ is a choice of moment map. A change
of~$\mum_0$ just adds a constant to~$t$ (i.e.\ changes~$L'$ by  
a homothety) and does not change the isometry class of
$\omega_\varphi$. 

In adapted coordinates, 
\begin{equation} \label{adapted}
\omega_\varphi = \frac{\sqrt{-1}}{2}
\left({{\varphi(\mum)}\over{2z^0\bar{z}^0}}\,dz^0\wedge d\bar{z}^0 
+\Bigl[g_M(\mum)\Bigr]_{\alpha\bar\beta}\,
dz^\alpha \wedge d\bar{z}^\beta\right)
\end{equation} 
along a fibre. In terms of globally defined functions, the metric
splits {\em along a fibre} into
\begin{equation}
\label{eqn:splitting}
\omega_{\rm fibre}+\omega_{\rm horiz}=ds\wedge d^cs+\omega_M(\mum)
\end{equation}
The expression~(\ref{fibg}) for the fibre metric follows at once, as
does the fact that~$s$ is the geodesic distance in the fibre.

\subsubsection*{The Ricci form}

From (\ref{adapted}), the volume form is
\begin{equation}\label{eqn:volume}
{{\omega_\varphi^{m+1}}\over{(m+1)!}}=(\varphi Q)(\mum)\,
\det\Bigl[(g_M)_{\alpha\bar\beta}\Bigr]\frac{1}{z^0\bar{z}^0}
\left(\Bigl({{\sqrt{-1}\over2}}\Bigr)^{m+1}
\prod_{i=0}^m dz^i\wedge d\bar{z}^i\right),
\end{equation}
from which (\ref{volumef}) follows immediately.

Denote by $\del_M$ and $\delbar_M$ the $\del$ and $\delbar$-operators
on $M$; it is necessary here to distinguish them from the
corresponding operators on $L$, which we continue to denote by $\del$
and $\delbar$. For every smooth function~$u$ on $I\times M$, 
$$
dd^c u = d_M d_M^c u 
- \varphi\frac{\del u}{\del \mum}\,\gamma
+ \frac{1}{\varphi}\frac{\del}{\del\mum}
\left[\varphi\frac{\del u}{\del\mum}\right]
d\mum\wedge d^c\mum + \mbox{cross-terms}. 
$$
(The cross-terms take the form $\delbar\mum\wedge \del_M (\del
u/\del\mum) + \mbox{ complex conjugate}$.) In order to take the trace
with respect to $\omega_\varphi$, either work locally or multiply by
$\omega_\varphi^m$. For the latter,
$$
\omega_\varphi^m = \omega(\mum)^m + m\,\omega(\mum)^{m-1}
\cdot {1\over\varphi}\,d\mum\wedge d^c\mum.$$
Wedging with $dd^c u$, the cross-terms drop out, so upon division
by~$\omega_\varphi^{m+1}$,
$$
\lapl_\varphi u = \lapl_{\omega_M(\mum)} u(\mum,\cdot) 
+ \frac{1}{2}\left(\frac{\del}{\del\mum}
\Bigl(\varphi\frac{\del u}{\del\mum}\Bigr) 
- \Bigl(\tr_{\omega_M(\mum)}\gamma\Bigr)\, 
\varphi \frac{\del u}{\del\mum}\right).
$$
The right-hand side of this is simplified by the following
observation: If the eigenvalues of~$\Beta$ are denoted $\beta_\nu(z)$,
$\nu=1,\ldots,m$, then
$Q(\mum)=\prod\limits_\nu\bigl(1-\mum\,\beta_\nu(z)\bigr)$, so
\begin{equation}
\label{eqn:dlogQ}
\tr_{\omega_M(\mum)}\gamma
= \tr\bigl[({\rm I}-\mum\Beta)^{-1}\Beta\bigr]
=\sum_{\nu=1}^m {{\beta_\nu(z)}\over{1-\mum\,\beta_\nu(z)}}
=-{\del\over{\del\mum}}(\log Q).
\end{equation}
Equations (\ref{boxphi})~and (\ref{sigma2}) follow immediately,
while~(\ref{sigma3}) follows from the observation
$$
\lapl_{\omega_M(\mum)} \log\varphi Q
= \lapl_{\omega_M(\mum)}(\log\varphi + \log Q) =
\lapl_{\omega_M(\mum)} \log Q
$$
($\varphi$ being pulled back by the first projection $I\times M \to
I$) and the formulae
\begin{eqnarray}
-\sqrt{-1}\del_M\delbar_M \log Q &=& \rho_{\omega_M(\mum)}-\rho_M \\
-\lapl_{\omega_M(\mum)} \log Q
&=&\tr_{\omega_M(\mum)}\bigl(\rho_{\omega_M(\mum)} -\rho_M\bigr)
\ =\ \sigma_M(\mum)-R(\mum).
\label{eqn:traces}
\end{eqnarray}

\subsubsection*{Total geodesy of fibres}

It is well-known that the fibres of~$L'$ are totally geodesic with
respect to~$\omega$. The proof given here serves as an excuse to
calculate the Levi-Civita connection of~$g$.

\begin{proposition}
\label{prop:tg}
Let $g_\varphi$ be a bundle-adapted metric on $L'\subset L$. Then each
fibre of~$L'$ is totally geodesic with respect to~$g_\varphi$. 
\end{proposition}

\begin{proof}  It suffices to calculate the Levi-Civita connection $D$
of $g_\varphi$ and show that $D_{\del_0}\del_0$ is tangent to the
fibre.  Let $(z^0,z)$ be a line bundle chart, $w^0=\log z^0$, and let
$\del_\alpha$ denote partial differentiation. Recall that the
connection form of $(L,h)$ is equal to $\Theta=\del\log r$. With
respect to the coordinates $(w^0,z)$, the vector-valued $(1,0)$-form
$\Theta$ is given by the column 
$$[\Theta_i]^{\rm t}=[\matrix{1 & \Theta_\alpha \cr}]^{\rm t},\qquad
\Theta_\alpha=h^{-1}\del_\alpha h,$$
and $\delbar(\Theta_\alpha\,dz^\alpha)=p^*\gamma$. The components
of~$g_\varphi$ are given by the Hermitian $(1+m)\times(1+m)$ block
matrix 
\begin{equation}
\label{eqn:tg1}
G=\left[\matrix{g_{0\bar0} & g_{0\bar\beta} \cr &\cr
g_{\alpha0} & g_{\alpha\bar\beta}\cr}\right]
=2\varphi(\mum)\left[\matrix{1 & [\Theta_{\bar\beta}]^{\rm t} \cr &\cr
[\Theta_\alpha] & [\Theta_\alpha][\Theta_{\bar\beta}]^{\rm t} \cr}
\right] + \left[{\matrix{0 & 0 \cr   &\cr 
0 & [g_M(\mum)]_{\alpha\bar\beta}\cr}}\right].
\end{equation}
The inverse matrix $G^{-1}$ is found by (block) row-reduction:
\begin{equation}
\label{eqn:tg2}
G^{-1}={1\over{2\varphi(\mum)}}\left[
\matrix{1 & \quad 0 \cr  &\cr 0  & \quad 0 \cr} \right] +
\left[\matrix{[\Theta_{\bar\beta}]^{\rm t}\,
[g_M(\mum)]^{\bar{\beta}\alpha}\, [\Theta_\alpha] &
-[\Theta_{\bar\beta}]^{\rm t}\,[g_M(\mum)]^{\bar{\beta}\alpha}
\cr &\cr
-[g_M(\mum)]^{\bar{\beta}\alpha}\, [\Theta_\alpha] &
\Neg[g_M(\mum)]^{\bar{\beta}\alpha} \cr}\right].
\end{equation}
The Levi-Civita connection form of~$g$ is represented in a line bundle
chart by the matrix-valued $(1,0)$-form $G^{-1}\del G$. A short
calculation shows that in adapted coordinates,
$$G^{-1}\del G=\left[\matrix{
(1/2)\varphi'(\mum)\del t & \del[\Theta_{\bar\beta}]^{\rm t} \cr & \cr 
\varphi(\mum)[g_M(\mum)]^{{\bar\beta}\alpha}\del[\Theta_\alpha] 
& [g_M(\mum)]^{{\bar\beta}\nu}
[g_M(\mum)]_{\alpha\bar\beta}\cr}\right]$$ along the
fibre. Evaluating this matrix-valued $(1,0)$-form on the tangent
vector $\del_0=\del/\del w^0$ at a point of the fibre gives the
representation (with respect to the frame $\{\del_i\}_{i=0}^n$) of the
covariant derivative $D_{\del_0}$ in the fibre direction. The
resulting matrix has first column $\left[\matrix{(1/2)\varphi'(\mum) &
0 & \cdots & 0\cr}\right]^{\rm t}$, which implies the covariant
derivative $D_{\del_0}\del_0$ is tangent to the fibre, i.e.\ that the
fibre is totally geodesic.
\end{proof}

\subsection{Completeness and Extendability of Fibre Metrics}
\label{section:complete}

In this section, the completeness properties of the metric~$g_\varphi$
are given in terms of the boundary behaviour of~$\varphi$.  Because
each fibre is totally geodesic by Proposition~\ref{prop:tg},
$g_\varphi$ is complete iff
\begin{itemize}
\item The metric $g_M(\mum)$ is complete for every $\mum\in I$, and
\item The fibre metric $g_{\rm fibre}$ is complete. 
\end{itemize}
The first condition will be assumed, so the task is to relate the
second to the boundary behaviour of the profile. By
equation~(\ref{eqn:mu}), completeness is guaranteed by divergence of
the $s$~integral, though divergence is not necessary since a
fibre metric may extend smoothly to the origin. In any case, criteria
in terms of order of vanishing are easier to work with.  By
definition, a profile is smooth on the closure of~$I$, so at each
finite endpoint~$\alpha$ there is an expansion
\begin{equation} \label{frob}
\varphi(\mum) = a_\ell(\mum -\alpha)^\ell 
+ O\bigl((\mum-\alpha)^{\ell+1}\bigr),\qquad a_\ell\not=0,
\end{equation}
with $\ell\geq0$ an integer.  For convenience, it will also be assumed
that~$\varphi$ is asymptotic to an integer power of~$\mum$ as
$\mum\to\infty$. This is no real loss, since the profiles of greatest
interest are {\em rational}\/ functions, which arise when
$g_\varphi$ has constant scalar curvature (or is formally
extremal). For such profiles, the relation between boundary behaviour
and completeness is very simple:
\begin{proposition}
\label{prop:extend}
Let $\varphi:I\to\R$ be a profile as above. Then the associated fibre
metric is complete if and only if one of the following conditions
holds at each endpoint of~$I$:
\begin{itemize}
\item Finite Endpoint{\rm(}s{\rm)}
\begin{description}
\item{\rm(i)} The profile $\varphi$~vanishes to first order, and
$|\varphi'|=2$ at the endpoint; or
\item{\rm(ii)} The profile vanishes to order at least two.
\end{description}
\item Infinite Endpoint{\rm(}s{\rm)}
\begin{description}
\item{\rm(iii)} The profile grows at most quadratically, i.e.\ there
is a positive constant~$K$ such that $\varphi(\mum)\leq K\mum^2$ for
$|\mum|\gg0$.
\end{description}
\end{itemize}
\end{proposition}

Apart from the proof of this result, this section is devoted to the
classification of $S^1$-invariant metrics of constant Gaussian
(equivalently scalar) curvature on $S^1$-invariant domains
of~$\P^1$. The two subsections may be read in either order; it
is hoped that Table~\ref{table-big} below will illuminate the proof of
Proposition~\ref{prop:extend}.

\subsubsection*{Proof of Proposition~\protect{\ref{prop:extend}}}

Lemma~\ref{lemma:interval} below asserts there is no loss of
generality in taking the lower endpoint of~$I$ to be~$0$; this is
assumed here for convenience.  The boundary condition~$\varphi(0)=0$
is necessary for completeness of~$g_{\rm fibre}$; if $\varphi(0)>0$
then the $s$~integral in~(\ref{eqn:mu}) converges at $\mum=0$, so the
distance to the level set $\{\mum=0\}$ is finite, while the length of
the corresponding $S^1$-orbit is $2\pi\sqrt{\varphi(0)}>0$. This is
impossible for a complete metric.

Hence $\varphi$ vanishes to order $\ell\geq1$, the $t$~integral
diverges at $\mum=0$, and the level set $\{\mum=0\}$ intersects the
fibre at the origin. If $\ell\geq 2$ the $s$~integral is also
divergent, so the origin is at infinite distance, i.e.\ the
fibre metric is complete. If $\ell=1$, the $s$-integral is convergent
and more work is needed to determine whether the fibre metric extends
smoothly to the level set $\{r=0\}$.

In a line bundle chart, $r =h\,|z^0|^2$, and by~(\ref{fibg}) the
fibre metric is
$$
g_{\rm{fibre}} = \varphi(\mum)\,\left|\frac{dz^0}{z^0}\right|^2 
= \left[{{\varphi(\mum)}\over r}\right]\,h\,\big|dz^0\big|^2.
$$
Since the horizontal part is smooth, $g_\varphi$ is smooth at $r=0$
iff $\varphi(\mum)/r$ has a finite, positive limit as $\mum \to0$.
(This limit makes sense, for~(\ref{eqn:mu}) gives~$t$ and hence~$r$ as
a function of $\mum$.)  Using (\ref{eqn:mu}) to find $t$, with
$\varphi(\mum)=a_1\mum+O(\mum^2)$, gives
$\log\mum ={\rm constant}+ a_1 t + \cdots$, or 
$\mum = a\,r^{a_1/2} + \cdots$ for some $a>0$. Thus $\varphi(\mum)/r$
has a finite, positive limit as $r\to0$ iff $a_1 =2$, i.e.\
$\varphi(0) =0$, $\varphi'(0)=2$.  This completes the analysis at the
lower endpoint of $I$. The analysis at the upper endpoint~$b$ is
similar when $b<\infty$; in particular, the fibre metric extends to
infinity in the fibre iff $\varphi(b)=0$ and $\varphi'(b)=-2$.

It remains to consider the case $(0,\infty)\subset I$. Completeness of
the fibre metric, i.e.\ unboundedness of the distance to the level set
$\{t=t_2\}$, limits the growth of $\varphi$ at~$\infty$ via the
relation
\begin{equation}
\label{eqn:t-complete}
\lim_{\mum\to\infty}\int^\mum
{{dx}\over{\sqrt{\varphi(x)}}}=\infty.
\end{equation}
The quadratic growth condition~(iii) is immediate; this completes the
proof of Proposition~\ref{prop:extend}. \QED

The boundary behaviour of $\varphi$ and  completeness properties of
$g_{\rm fibre}$ are summarized in Table~\ref{table-complete}.  
\begin{table}[hbt]
\begin{centering}
\begin{tabular}{|l|c|c|c|c|} 
\hline
Profile Near $\mum=0$ 
		& Distance to $\mum=0$ & End Geometry & Complete/Smooth \\
\hline
\quad$\vphantom{\Big|}\varphi(\mum)=2\mum+O(\mum^2)$ 
			& Finite & Smooth Extension & Yes\\
\quad$\vphantom{\big|}0<\varphi'(0)\neq2$ & Finite 
			& Point Singularity & No\\
\quad$\vphantom{\Big|}\varphi(\mum)\leq K\mum^2$ & Infinite
			& Finite-area Cusp & Yes\\
\hline\hline
Profile Near $\mum=\infty$ 
		& Distance to $\mum=\infty$ & End Geometry & Complete \\
\hline
\quad$\vphantom{\Big|}\varphi(\mum)\to0$ 
			& Infinite & Infinite-area Cusp & Yes \\
\quad$\varphi(\mum)\to{\rm const}>0$ & Infinite	& Cylindrical & Yes \\
\quad$\vphantom{\Big|}\varphi(\mum)\sim K\mum$ & Infinite 
					& Planar/Conical & Yes \\
\quad$\varphi(\mum)\sim K\mum^2$ & Infinite & Hyperbolic & Yes \\
\quad$\vphantom{\Big|}\varphi(\mum)\geq K\mum^3$ & Finite &  & No \\
\hline
\end{tabular}
\caption{Completeness and boundary/asymptotic behavior for rational
profiles.}\label{table-complete}
\end{centering}
\end{table}

\subsubsection*{Geometry of fibre metrics} 

In this section, the classification of complete $S^1$-invariant
metrics on subsets of $\P^1$ is considered from the point of view of
the momentum construction.  The results are summarized in
Table~\ref{table-big} below, which shows clearly how much simpler
these metrics look in `momentum coordinates' than in the standard
complex coordinate~$z$.

%% Table~\ref{table-big}.  
\begin{table}[hbt]
\label{page-big}
\begin{centering}
\begin{tabular}{||c|c|c|c|c|c|c|c||} 
\hline
  & & & & \small Distance  to &
{\small Area  near} & \small Metric & 
 \\
& $I^0$ & $\varphi(\mum)$ & $r$-range & 
$\mum = a$ & $\mum = a$ & $\sigma$ & \small Domain\\
& & & & $\mum = b$ & $\mum = b$ & 
$\vphantom{\Big|}{{\varphi(\mum)}\over r}$ & \\
\hline
 & & & & & &  \small Fubini-Study &  \\
(i) & $(0,2c^{-2})$ & $2\mum - c^2\mum^2$ & $[0,\infty]$ & finite &
finite & $c^2$ & $\P^1$ \\
 & & & & finite & finite & $\vphantom{\Big|}\frac{4}{c^2(1+r)^2}$ & \\
\hline
& & & & & & \small flat plane & \\
(ii) & $(0,\infty)$ & $2\mum$ & $[0,\infty)$ & finite &
finite & $0$ &  $\C$ \\  & & & & infinite & infinite & 
$\vphantom{\Big|}1$ & \\
\hline
 & & & & & &  \small Poincar\'e disc&  \\
(iii) & $(0,\infty)$ & $2\mum + c^2\mum^2$ & $[0, 1)$ & finite &
finite & $- c^2$ & $\Delta$ \\  & & & & infinite & infinite 
	& $\vphantom{\Big|}\frac{4}{c^2(1-r)^2}$ & \\
\hline
& & & & & & \small hyperbolic cusp &  \\
(iv) & $(0,\infty)$ & $c^2\mum^2$ & $(0, 1)$ & infinite &
finite & $- c^2$ & $\Delta^\times$ \\
 & & & & infinite & infinite &  
$\vphantom{\Big|}\frac{4}{c^2\,r(\log r)^2}$ & \\
\hline
& & & & &  &  \small flat cylinder &  \\ 
(v) & $\R$ & $\alpha^2$ & $(0, \infty)$ & infinite &
infinite & $0$ & $\C^\times$ \\
 & & & & infinite & infinite & $\vphantom{\Big|}\frac{\alpha^2}{r}$ & \\
\hline
& & & & & & \small hyperbolic annulus &  \\
(vi) & $\R$ & $\alpha^2 + c^2\mum^2$ &
$(e^{-{{\pi}\over{c\alpha}}}, e^{{\pi}\over{c\alpha}})
$ & infinite & infinite & $- c^2$ & \small Annulus \\ 
 & & & & infinite & infinite & 
$\vphantom{\Big|}\frac{\alpha^2}{r\cos^2(c\alpha\log r/2)}$ & \\ 
\hline
\end{tabular}
\caption{Classification of complete, circle-invariant metrics with
constant scalar curvature on domains in $\P^1$ by momentum
data.}\label{table-big} 
\end{centering}
\end{table}

Since the scalar curvature is given in terms of the profile as
$-\varphi''/2$, constancy of the scalar curvature implies that the
profile is quadratic. The list of profiles, together with the
intervals $I^0$ on which they are positive, appear in the first two
columns of the table. As usual, the normalization $\inf I = 0$ has
been used in all but cases (v) and (vi), where $I =\R$. Note that $I^0$
is the {\em interior} of the image~$I$ of the moment map.

Performing the $t$~integral in each case yields $t$---and hence
$r=|z|^2$---as a function of~$\mum$. From these formulae, the
$r$~interval in the third column follows at once. The next three
columns follow either directly or from Proposition~\ref{prop:extend}.
Finally, the metric in question is named, its scalar curvature is
given, and the conformal factor $[\varphi(\mum)/r]$ is written
explicitly as a function of $r$.  The latter is obtained by inverting
the dependence of $r$ upon $\mum$ and substituting in the formula for
$\varphi$.

\begin{remark} (iv) is the finite-area hyperbolic cusp. In (vi) the
annulus is determined, up to conformal equivalence, by the ratio of
the inner and outer radii and hence by 
$c\alpha$. In the examples (i)--(iii), with
$\varphi'(0)=2$, metrics with conical singularity at the origin arise
by changing the coefficient of $\mum$ to $a\not=2$.
\end{remark}

It is instructive to compare these metrics with surfaces of
revolution in~$\R^3$.  Let $\xi$ be a positive function, whose graph
sweeps out a surface of revolution. The area and length elements are
given in terms of the independent variable~$y$ as
$$d\mum = \xi(y)\sqrt{1+\xi'(y)^2}\,dy,\qquad
ds = \sqrt{1+\xi'(y)^2}\,dy.$$ 
The profile is the length squared of the vector field~$X$, or
$\varphi(\mum)=\xi(y)^2$. Differentiating and using the previous
equations gives
$$\varphi'(\mum)={{2\xi'(y)}\over{\sqrt{1+\xi'(y)^2}}},\qquad
{\rm or}\qquad 
\xi'(y)={{\varphi'(\mum)}\over{\sqrt{4-\varphi'(\mum)^2}}},$$ 
which implies $|\varphi'(\mum)|\leq 2$, with equality iff
$|\xi'(y)|=\infty$: The fibre metric embeds as a surface of revolution
iff the profile is not too steep.  Comparing with the profiles in the
table, one notes the well-known facts that the portion of the
cusp~(iv) corresponding to the momentum sub-interval $(0,1/c^2)$ 
embeds as a surface of revolution in $\R^3$ (the pseudosphere), while
no annulus in the Poincar\'e disc~(iii) arises in this way.

It is also worth noting that intuition deriving from surfaces of
revolution can be misleading. Taking $I=[0,\infty)$ and
$\varphi(\mum)=2\mum+\mum^3$ gives an incomplete metric of
infinite area on the disk!

\section{Metrics of Infinite Volume}
\setcounter{equation}{0}
\setcounter{figure}{0}
\label{section:inf-vol}

This section is devoted to the proofs of Theorems~\ref{thm:A},
\ref{thm:B}, and \ref{thm:C}. Section~\ref{section:line} deals with
Theorems~\ref{thm:A} and \ref{thm:B}, while
Section~\ref{section:vbundle} with Theorem~\ref{thm:C}. The proofs are
separated in this way because of the behaviour of the family of
K\"ahler forms~$\omega_M(\tau)$; in Theorems~\ref{thm:A} and
\ref{thm:B}, $\omega_M(0)$ is non-degenerate, whereas in
Theorem~\ref{thm:C} the family~$\{\omega_M(\mum)\}$ drops rank at
$\mum=0$, leading to a metric on a partial blow-down of~$L$.

\subsection{Metrics on Line Bundles}
\label{section:line}

Let $\{p:(L,h)\to(M,\omega_M),\ I\}$ be $\sigma$-constant horizontal
data with $(L,h)$ not flat, and with $I= [0,\infty)$ or
$(0,\infty)$. In particular, $\gamma\leq0$, so
$\omega_M(\mum)=\omega_M-\mum\gamma$ is a K\"ahler form on~$M$ for all
$\mum\geq0$.  Further, the metric $g_M(\tau)$ is complete by hypothesis
because $g_M$ is assumed to be complete and $-\tau\gamma$ is positive
semidefinite if $\tau \in I$.  As in Section~\ref{section:KS}, define
functions $Q$,~$R:I\to\R$ by
$$
Q(\mum) = \det (1-\mum\Beta),\qquad 
R(\mum) = \tr\Bigl[(1 -\mum\Beta)^{-1}\varrho\Bigr].
$$
Because the data are $\sigma$-constant, $Q(\mum)$ is the (constant)
scale factor for the volume form of $\omega_M(\mum)$ relative
to~$\omega_M$ and is a polynomial in~$\mum$ with only negative roots,
while $R(\mum)=\sigma_M(\mum)$ is the scalar curvature
of~$g_M(\mum)$ and is a rational function that is bounded below
on~$I$. Indeed, $R$~has an asymptotic value $R(\infty)$ as
$\mum\to\infty$, which may be interpreted as the trace of $\varrho$
restricted to the $0$-eigenbundle of~$\Beta$.  By 
Corollary~\ref{cor:CKS}, the scalar curvature of~$\omega_\varphi$ is
\begin{equation} \label{pres}
\sigma_\varphi=R(\mum)-(1/2Q)\bigl(\varphi Q\bigr)''(\mum).
\end{equation}
Now let $\sigma$ be a function on~$I$. The problem of prescribing the
scalar curvature of $\omega_\varphi$ is given by the equation $\sigma
= \sigma_\varphi$, which has the solution
\begin{equation}
\label{eqn:gen-profile}
(\varphi Q)(\mum) = (\varphi Q)(0)+(\varphi Q)'(0)\mum 
+ 2\int_0^\mum (\mum-x)\Bigl(R(x)-S(x)\Bigr)Q(x)\,dx
\end{equation}
in terms of the initial data $\varphi(0)$~and $\varphi'(0)$. The
momentum construction yields a metric of infinite fibre area iff
$\varphi$ is positive on $(0,\infty)$. This metric is complete if
$\varphi$ grows at most quadratically at~$\infty$ and satisfies the
boundary conditions given in Proposition~\ref{prop:extend}:
$\varphi(0) = 0$, and either $\varphi'(0)=2$ or $\varphi'(0) =0$. The
first case is the one needed for Theorem~\ref{thm:A}, the second for
Theorem~\ref{thm:B}. These will now be considered in turn.

\subsubsection*{Proof of Theorem \protect\ref{thm:A}}

Setting $\sigma = c$ (constant) and using the initial conditions
$\varphi(0) =0$, $\varphi'(0)=2$ in (\ref{eqn:gen-profile}) gives 
\begin{equation}
\label{eqn:varphi}
\varphi(\mum) = {2\over{Q(\mum)}}\left(
\mum + \int_0^\mum (\mum-x)\Bigl(R(x)-c\Bigr)Q(x)\,dx\right).
\end{equation}
The notation $\varphi_c(\mum)$ or $\varphi(\mum,c)$ will be used when
the dependence of the profile on~$c$ is being emphasized.  Define the
set ${\rm J}\subset\R$ of ``allowable scalar curvatures'' by
$${\rm J} = \{c\in\R\mid\varphi_c(\mum) > 0 \mbox{ for all } \mum>0\}.$$
In words, ${\rm J}$ is the set of~$c$ for which
equation~(\ref{eqn:varphi}) defines a momentum profile on~$I$.

\begin{lemma}
\label{lemma:J}
There is a $c_0\in\R$ such that either ${\rm J}=(-\infty,c_0)$ or
${\rm J}=(-\infty,c_0]$.
\end{lemma}

\begin{proof}
The function $R$ is bounded on~$I$, so the integrand in
(\ref{eqn:varphi}) is positive for $c\ll0$, implying $\varphi_c>0$
on~$(0,\infty)$ for $c\ll0$. In particular, ${\rm J}$~is non-empty.
If $\mum>0$, then
$${{\del\varphi}\over{\del c}}(\mum,c)=-{2\over{Q(\mum)}}
\int_0^\mum (\mum-x) Q(x)\,dx < 0,$$
i.e.\ $\varphi(\mum,c)$ is strictly decreasing with respect to~$c$.
Consequently, if $c\in {\rm J}$~and $c'<c$, then $c'\in {\rm
J}$. Finally, ${\rm J}$ is bounded above since by~(\ref{eqn:varphi}),
$\varphi_c$ is not everywhere positive for $c>R(\infty)$. In summary,
${\rm J}$ is a half-line, unbounded below.  Set $c_0=\sup {\rm J}\leq
R(\infty)$.
\end{proof}

Since the initial condition ensures that $\varphi_c$ is positive for
sufficiently small positive $\tau$, the momentum construction yields a
metric which for notational convenience will be denoted~$g_c$. 

\begin{lemma}
\label{lemma:complete} If $c < c_0$ then $g_c$ is a complete
metric on $\Delta(L)$, and $g_{c_0}$ is a complete metric on~$L$.
\end{lemma}

\begin{proof}  
By Proposition~\ref{prop:extend} and Table~\ref{table-complete},
completeness of the metric is equivalent to `at most quadratic' growth
of~$\varphi_c$ at~$\infty$, provided the profile is positive
on~$(0,\infty)$. To establish the growth condition, it is easiest to
write~$R-c$ as a constant plus a rational function~$R_0$ vanishing
at~$\infty$:
$$R(\mum)-R(\infty)=:R_0(\mum),\quad c-R(\infty)=:\tilde{c}, 
\qquad\mbox{so}\quad R(\mum)-c=R_0(\mum)-\tilde{c}.$$
Define polynomials
$$P_1(\mum) =  \mum + \int_0^\mum (\mum-x) R_0(x) Q(x)\,dx,\qquad
  P_2(\mum) =  \int_0^\mum (\mum-x) Q(x)\,dx;$$
Because $R_0$ vanishes at~$\infty$, $\deg P_1\leq 1+\deg Q$, and $\deg
P_2=2+\deg Q$. From (\ref{eqn:varphi}),
\begin{equation}
\label{eqn:c-dep}
\varphi(\mum,c) 
={2\over{Q(\mum)}}\Bigl(P_1(\mum)-\tilde{c} P_2(\mum)\Bigr),
\end{equation}
so if $\tilde{c}=0$, then $\varphi_c(\mum)\leq K\mum$ as
$\mum\to\infty$, while $\varphi_c(\mum)\sim K\mum^2$ if $\tilde{c}<0$,
i.e.\ if $c<R(\infty)$. Since $c_0\leq R(\infty)$, it follows that
each profile~$\varphi_c$ with $c<c_0$ gives rise to a complete fibre
metric, hence to a complete metric~$g_c$. Because $\varphi_c$~grows
quadratically, the $t$~integral converges as $\mum\to\infty$, so up to
homothety $g_c$~lives on the unit disk bundle~$\Delta(L)$.

It remains to investigate the borderline case. First observe that
every profile~$\varphi_c$ with $c<c_0$ is bounded away from zero
except near $\mum=0$ (since $\varphi_c'(0)=2$, and~$\varphi_c$ is
positive on~$(0,\infty)$ and has infinite limit as $\mum\to\infty$).
By the proof of Lemma~\ref{lemma:J}, $\varphi(\mum,c)$ is decreasing
in~$c$.  Since $c_0$ is the supremum of~$c$ for which $\varphi_c$ is
positive on~$(0,\infty)$, it follows that~$\varphi_{c_0}$ is
non-negative on~$I$ by continuity of~$\varphi(\mum,c)$ in~$c$.  The
borderline profile is not identically zero, since
$\varphi_{c_0}'(0)=2$.

Two possibilities occur: $\varphi_{c_0}$ has a positive zero, or
is positive on $(0,\infty)$. In the first case let $b$ be the first
positive zero of $\varphi_{c_0}$. Then $\varphi_{c_0}'(b)=0$ as well,
for $\varphi_{c_0}$ is real-analytic and non-negative. Since
$\varphi_{c_0}$ vanishes to order at least two, the $t$~and
$s$~integrals diverge, so the associated metric lives on the total
space of~$L$ and is complete, but has finite-area fibres, see also
Table~\ref{table-complete}.

Consider now the second possibility, $\varphi_{c_0} >0$ on
$(0,\infty)$. We claim in this case that $c_0 = R(\infty)$ (i.e.\
$\tilde{c}=0$). For if not, the borderline profile is positive and
grows quadratically, hence is bounded away from zero except near
$\mum=0$; a glance at~(\ref{eqn:c-dep}) shows $\tilde{c}$ may be
increased slightly, preserving positivity of the profile, but this
contradicts the definition of~$c_0$. Hence
$\varphi_{c_0}(\tau)=(P_1/Q)(\mum) \leq K\tau$ and again $g_{c_0}$
lives on~$L$. 
\end{proof}

\begin{figure}[htb]  %% C.pic
\setlength{\unitlength}{0.00083333in}
\begingroup\makeatletter\ifx\SetFigFont\undefined
% extract first six characters in \fmtname
\def\x#1#2#3#4#5#6#7\relax{\def\x{#1#2#3#4#5#6}}%
\expandafter\x\fmtname xxxxxx\relax \def\y{splain}%
\ifx\x\y   % LaTeX or SliTeX?
\gdef\SetFigFont#1#2#3{%
  \ifnum #1<17\tiny\else \ifnum #1<20\small\else
  \ifnum #1<24\normalsize\else \ifnum #1<29\large\else
  \ifnum #1<34\Large\else \ifnum #1<41\LARGE\else
     \huge\fi\fi\fi\fi\fi\fi
  \csname #3\endcsname}%
\else
\gdef\SetFigFont#1#2#3{\begingroup
  \count@#1\relax \ifnum 25<\count@\count@25\fi
  \def\x{\endgroup\@setsize\SetFigFont{#2pt}}%
  \expandafter\x
    \csname \romannumeral\the\count@ pt\expandafter\endcsname
    \csname @\romannumeral\the\count@ pt\endcsname
  \csname #3\endcsname}%
\fi
\fi\endgroup
{ %%\renewcommand{\dashlinestretch}{30}
\begin{picture}(7362,3177)(0,-10)
%% Right-hand graph labels
%% Origin of right-hand graph at (4350,450)
\put(4200,3000){\makebox(0,0)[lb]{\smash{{{$c$}}}}}
\put(4125,1300){\makebox(0,0)[lb]{\smash{{{$c_0$}}}}}
\put(5850,1725){\makebox(0,0)[lb]{\smash{{{$c=C(\mum)$}}}}}
\put(4425,775){\makebox(0,0)[lb]{\smash{{{\SetFigFont{12}{14.4}{rm}J}}}}}
\put(4325,195){\makebox(0,0)[lb]{\smash{{{$0$}}}}}
\put(2700,195){\makebox(0,0)[lb]{\smash{{{$\mum$}}}}}
%% Left-hand graph
\path(225,3150)(225,3148)(226,3144)
	(226,3137)(228,3126)(229,3111)
	(231,3091)(234,3067)(237,3039)
	(240,3007)(244,2971)(248,2934)
	(253,2895)(257,2855)(261,2816)
	(265,2776)(270,2739)(274,2702)
	(277,2667)(281,2635)(285,2604)
	(288,2575)(291,2548)(294,2523)
	(297,2498)(300,2475)(303,2448)
	(307,2422)(310,2395)(314,2370)
	(317,2344)(320,2319)(324,2294)
	(327,2269)(331,2244)(334,2219)
	(338,2195)(341,2171)(344,2147)
	(348,2123)(351,2099)(355,2076)
	(358,2054)(361,2032)(365,2010)
	(368,1989)(372,1969)(375,1950)
	(380,1925)(384,1902)(389,1880)
	(393,1859)(397,1840)(401,1821)
	(405,1803)(409,1786)(413,1769)
	(417,1752)(422,1735)(426,1718)
	(431,1701)(437,1684)(443,1667)
	(450,1650)(457,1634)(465,1618)
	(474,1603)(483,1587)(492,1572)
	(502,1557)(511,1542)(521,1527)
	(532,1512)(542,1498)(552,1484)
	(562,1471)(572,1458)(581,1446)
	(591,1435)(600,1425)(612,1413)
	(625,1404)(637,1396)(648,1389)
	(660,1382)(671,1377)(683,1372)
	(695,1367)(708,1362)(721,1358)
	(735,1354)(750,1350)(764,1348)
	(778,1346)(794,1345)(810,1344)
	(827,1343)(845,1343)(863,1343)
	(880,1343)(898,1343)(915,1344)
	(931,1345)(947,1346)(961,1348)
	(975,1350)(990,1354)(1004,1358)
	(1017,1363)(1028,1369)(1039,1374)
	(1050,1380)(1061,1387)(1072,1394)
	(1083,1401)(1096,1408)(1110,1416)
	(1125,1425)(1137,1432)(1149,1440)
	(1162,1448)(1176,1457)(1189,1466)
	(1203,1475)(1216,1484)(1230,1493)
	(1245,1502)(1259,1512)(1273,1522)
	(1288,1532)(1303,1542)(1318,1553)
	(1334,1564)(1350,1575)(1365,1586)
	(1380,1597)(1395,1608)(1410,1620)
	(1424,1632)(1439,1645)(1453,1658)
	(1468,1671)(1482,1684)(1497,1697)
	(1513,1710)(1529,1723)(1546,1737)
	(1564,1749)(1583,1762)(1604,1775)
	(1626,1787)(1650,1800)(1669,1809)
	(1690,1819)(1711,1828)(1733,1838)
	(1756,1848)(1780,1858)(1805,1868)
	(1830,1878)(1856,1888)(1882,1899)
	(1909,1909)(1936,1920)(1963,1930)
	(1990,1940)(2017,1950)(2044,1960)
	(2070,1969)(2097,1979)(2123,1987)
	(2149,1996)(2175,2004)(2200,2011)
	(2225,2018)(2250,2025)(2275,2031)
	(2301,2037)(2327,2042)(2355,2047)
	(2384,2051)(2415,2056)(2448,2060)
	(2482,2064)(2518,2068)(2555,2072)
	(2593,2076)(2630,2080)(2667,2084)
	(2701,2087)(2734,2090)(2763,2092)
	(2788,2095)(2809,2096)(2825,2098)
	(2836,2099)(2844,2099)(2848,2100)(2850,2100)
%% Left-hand graph labels
%% Origin of left-hand graph at (0,450)
\put(0,3000){\makebox(0,0)[lb]{\smash{{{$c$}}}}}
\put(-75,1300){\makebox(0,0)[lb]{\smash{{{$c_0$}}}}}
\put(7125,195){\makebox(0,0)[lb]{\smash{{{$\mum$}}}}}
\put(125,195){\makebox(0,0)[lb]{\smash{{{$0$}}}}}
\put(225,775){\makebox(0,0)[lb]{\smash{{{\SetFigFont{12}{14.4}{rm}J}}}}}
\put(1650,1500){\makebox(0,0)[lb]{\smash{{{$c=C(\mum)$,}}}}}
\put(1450,1300){\makebox(0,0)[lb]{\smash{{{i.e.\ $\varphi(\mum,c)=0$}}}}}
\thicklines
\path(4350,1350)(4350,450) %% Right-hand lower y-axis
\path(150,1275)(150,450)   %% Left-hand lower y-axis
\thinlines
\blacken\path(4380.000,3030.000)(4350.000,3150.000)(4320.000,3030.000)(4380.000,3030.000) 		  %% Right-hand y-axis arrow
\path(4350,3150)(4350,450) %% Right-hand upper y-axis
\blacken\path(180.000,3030.000)(150.000,3150.000)(120.000,3030.000)(180.000,3030.000)			  %% Left-hand y-axis arrow
\path(150,3150)(150,450)   %% Left-hand lower y-axis
%% Right-hand graph
\path(4425,3150)(4425,3148)(4425,3145)
	(4426,3138)(4426,3127)(4427,3111)
	(4428,3091)(4429,3065)(4430,3034)
	(4432,2998)(4434,2958)(4437,2913)
	(4440,2864)(4443,2813)(4446,2759)
	(4450,2704)(4454,2648)(4459,2592)
	(4463,2537)(4469,2484)(4474,2432)
	(4480,2382)(4487,2334)(4493,2289)
	(4501,2247)(4509,2208)(4518,2172)
	(4527,2138)(4538,2106)(4549,2077)
	(4562,2050)(4575,2025)(4592,1998)
	(4610,1972)(4629,1948)(4648,1925)
	(4669,1904)(4689,1884)(4710,1865)
	(4731,1847)(4752,1830)(4774,1814)
	(4795,1798)(4817,1783)(4839,1768)
	(4862,1754)(4885,1740)(4910,1727)
	(4935,1713)(4961,1700)(4988,1687)
	(5017,1674)(5047,1661)(5080,1648)
	(5114,1635)(5151,1622)(5190,1610)
	(5232,1598)(5277,1586)(5325,1575)
	(5357,1568)(5391,1562)(5427,1556)
	(5464,1550)(5504,1544)(5547,1539)
	(5591,1534)(5639,1528)(5689,1523)
	(5742,1518)(5797,1514)(5855,1509)
	(5916,1504)(5979,1500)(6044,1495)
	(6112,1490)(6181,1486)(6252,1482)
	(6324,1477)(6397,1473)(6471,1469)
	(6544,1465)(6616,1461)(6688,1457)
	(6757,1454)(6825,1450)(6890,1447)
	(6951,1444)(7009,1441)(7063,1438)
	(7112,1436)(7157,1434)(7196,1432)
	(7231,1430)(7261,1429)(7285,1428)
	(7306,1427)(7321,1426)(7333,1426)
	(7341,1425)(7346,1425)(7349,1425)(7350,1425)
\thicklines
\path(4350,1350)(4350,450)
\path(150,1275)(150,450)
\thinlines
\path(4350,450)(7350,450)
\blacken\path(7230.000,420.000)(7350.000,450.000)(7230.000,480.000)(7230.000,420.000)
\path(150,450)(2925,450)
\blacken\path(2805.000,420.000)(2925.000,450.000)(2805.000,480.000)(2805.000,420.000)
%% More left-hand graph labels
\put(4375,1350){\dashbox{60}(2975,0){}} %% (7350,1350)
\put(4315,1310){$\bullet$}
\put(4550,1425){$R(\infty)$}
\put(150,2175){\dashbox{60}(2725,0){}} %% (750,2175)
\put(110,1300){$\circ$}
\put(400,2250){$R(\infty)$}
\end{picture}}
\caption{The rational function~$C$ and
the interval~$\rm J$ of allowable scalar curvatures}
\label{fig:C}
\end{figure}
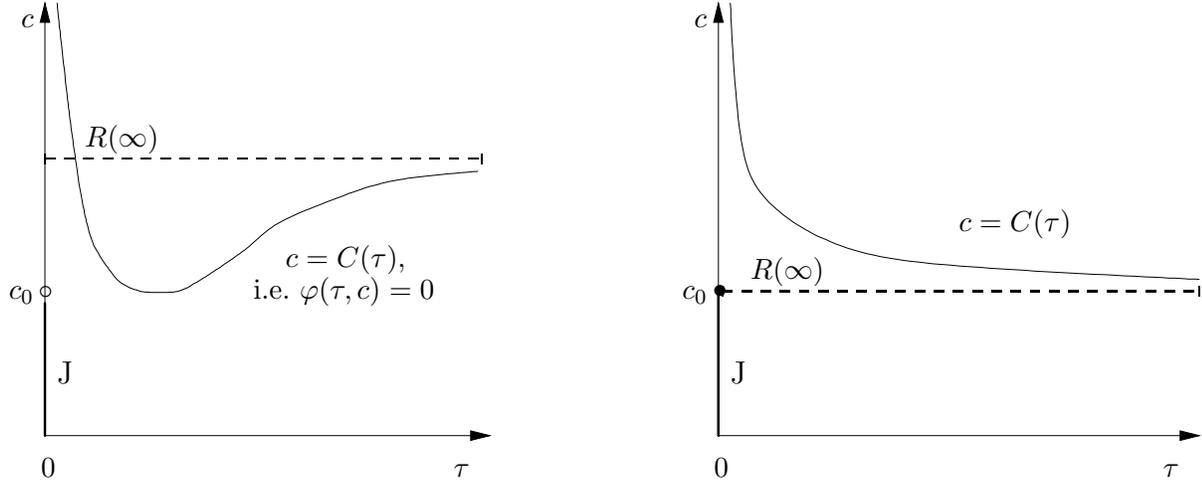

The dichotomy at $c=c_0$ is summarized as follows:
\begin{itemize}
\item $\varphi_{c_0}$ has a positive zero in~$I$ and yields a metric
  with fibrewise finite area. 
\item $\varphi_{c_0}$ is positive on $(0,\infty)$ and yields a
  complete metric of infinite fibre area; 
\end{itemize}
As shown above, the second alternative implies $c_0=R(\infty)$ (so
$c_0<R(\infty)$ implies the first alternative, see
Figure~\ref{fig:C}), but this is the only general conclusion that can
be drawn. Further, it is not necessarily true that $c_0=0$ (since
generally $R(\infty)\neq0$, for example); this issue is addressed in
detail below. However, if $\gamma<0$, or if the construction yields a
metric that is Einstein, then the borderline metric is scalar-flat,
while the others have negative curvature.

The borderline constant~$c_0$ has an alternative interpretation.
Consider the equation $\varphi(\mum,c)=0$. From (\ref{eqn:c-dep}),
this level set is the graph of the rational function $C:I\to\R$
defined by
$$C(\mum) = R(\infty) + \frac{P_1(\mum)}{P_2(\mum)}.$$
The degree of~$P_1$ is less than the degree of~$P_2$, so $C(\mum) \to
R(\infty)$ as $\mum\to\infty$. Furthermore, $P_1$~vanishes to order
one and $P_2$ vanishes to order two at $\mum=0$, so
$C(\mum)\to+\infty$ as $\mum\to 0^+$. Hence the function $C$~is
bounded below on $(0,\infty)$, and it follows immediately from the
definition that $c_0=\inf\{C(\mum)\mid \mum\in(0,\infty)\}$, see
Figure~\ref{fig:C}. 

If $c>c_0$, then $\varphi_c$ is {\em not}\/ non-negative on
$(0,\infty)$, so there is a first positive root~$b$, and the fibre
metric has finite area. Three possibilities occur:
\begin{enumerate}
\item $\varphi_c'(b) = 0$, and the fibre metric is complete, with a
  cusp end at $\mum=b$;
\item $\varphi_c'(b) = -2$, and the metric extends smoothly to the
  $\P^1$-bundle~$\widehat{L}$ (but note Remark~\ref{remark:positive}
  below); 
\item $\varphi_c'(b)\neq 0,\ -2$, so the metric is incomplete and has
  no smooth extension.
\end{enumerate}

\begin{lemma}
\label{lemma:finite}
There are at most finitely many $c>c_0$ for which the profile
$\varphi_c$ satisfies one of the first two conditions.
\end{lemma}

\begin{proof}
By~(\ref{eqn:c-dep}), $\varphi(\mum,c)=0$ if and only if
$\tilde{c}=P_1(\mum)/P_2(\mum)$. Substitution shows that if
$\varphi(\mum,c)=0$, then
$${{\partial\varphi}\over{\partial \mum}}(\mum,c)
  =2{{P_1(\mum)}\over{Q(\mum)}}
\Biggl(\log{{P_1}\over{P_2}}\Biggr)'(\mum).$$
This is a non-constant rational function, which takes the values
$0$~and $-2$ for at most finitely many values of~$\mum$. Consequently,
there are at most finitely many pairs $(\mum,c)$ satisfying
$$\varphi(\mum,c)=0\qquad{\rm and}\qquad
{{\del\varphi}\over{\del\mum}}(\mum,c)=0\ {\rm or}\ -2.$$
In any event, a value $c>c_0$ does not give rise to a metric with
infinite volume.
\end{proof}

The final task is to determine when a metric just constructed is
Einstein.  Equation~(\ref{ommu}) expresses the K\"ahler form
of~$g_\varphi$ in terms of $\varphi$~and the horizontal data,
while~(\ref{cricc}) similarly expresses the Ricci form. Setting
$\rho_\varphi = \lambda\,\omega_\varphi$ gives
$$
-\Bigl(\frac{1}{2Q}(\varphi Q)'\Bigr)' = \lambda
\qquad {\rm and}\qquad
\rho_M +\frac{1}{2Q}\bigl(\varphi Q\bigr)'(\mum)\,\gamma 
= \lambda(\omega_M - \tau\gamma).
$$
Integrating the first and using the initial conditions $\varphi(0)=0$,
$\varphi'(0)=2$, 
\begin{equation}
\label{eqn:einstein}
\frac{1}{2Q}(\varphi Q)'(\mum) = 1 - \lambda\tau.
\end{equation}
Substituting this back into the second equation,
$$
\rho_M + \gamma = \lambda\omega_M.
$$
Assume from now on that $\lambda\leq 0$, see
Remark~\ref{remark:positive} below. Integrating~(\ref{eqn:einstein}),
again  using $\varphi(0)=0$, gives 
\begin{equation}
\label{eqn:einstein-profile}
\varphi(\mum)={2\over{Q(\mum)}}\int_0^\mum (1-\lambda x)Q(x)\,dx,
\end{equation}
which is clearly positive for all $\mum>0$. Completeness follows since
$\varphi$ grows linearly (if $\lambda=0$) or quadratically (if
$\lambda<0$), and the scalar curvature is $\lambda(m+1)$.  Conversely,
a bundle-adapted metric arising in this way is Einstein:

\begin{lemma}
\label{lemma:corresp}
If the horizontal data are $\sigma$-constant and satisfy
$\rho_M+\gamma=\lambda \omega_M$ for some $\lambda\leq0$, then the
bundle-adapted metric with $c=\lambda(m+1)$ is Einstein-K\"ahler.
\end{lemma}

\begin{proof}  
Write $\rho_M=\lambda \omega_M-\gamma
=\lambda(\omega_M-\mum\gamma)-(1-\lambda\mum)\gamma$.  Taking the
trace with respect to $\omega_M(\mum)$, using the definition of~$R$
and recalling that $\tr_{\omega_M(\mum)}\gamma = -(\del/\del\mum)(\log
Q)$ by equation~(\ref{eqn:dlogQ}), immediately implies
$$R(\mum)= \lambda m + (1-\lambda \mum){{Q'(\mum)}\over{Q(\mum)}}.$$
If $c=\lambda(m+1)$, then
$$\Bigl(R(\mum)-c\Bigr)Q(\mum)
= Q'(\mum)-\lambda\Bigl(Q(\mum)+\mum Q'(\mum)\Bigr)
= {d\over{d\mum}}\Bigl((1-\lambda \mum)Q(\mum)\Bigr).$$ 
Integrating twice, using $\varphi(0)=0$ and $\varphi'(0)=2$, proves
the profile~(\ref{eqn:einstein-profile}) coincides with the
profile~(\ref{eqn:varphi}).  
\end{proof}

This completes the proof of Theorem~\ref{thm:A}.  

\begin{remark}
\label{remark:positive}
When the Einstein constant~$\lambda$ is positive, every completion of
$\lx$ is compact by Myer's Theorem, so there is a non-vacuous boundary
condition imposed on~(\ref{eqn:einstein}), namely that
$\varphi'(b)=-2$ when $\varphi(b)=0$. Koiso and Sakane have shown this
is equivalent to vanishing of the Futaki invariant of the
compactification. This boundary condition is never satisfied if the
curvature form~$\gamma$ is negative semi-definite unless
$\gamma=0$. To see this, observe that $\varphi(b)=0$ if and only if
$$\int_0^b(1-\lambda x)Q(x)\,dx=0,$$
in which case a short calculation shows that $\varphi'(b)=-2$ if and
only if $b\lambda=2$. Substituting back
into~(\ref{eqn:einstein-profile}), 
$$\varphi(\mum)={2\over{Q(\mum)}}\int_0^\mum (1-{2\over b}x)Q(x)\,dx 
= -{{4}\over{b\,Q(\mum)}}
\int_{-b/2}^{\mum-(b/2)}x\,Q\Bigl(x+(b/2)\Bigr)\,dx,$$ and
this is not zero when $\mum=b$ because $Q$ is positive and increasing
on the interval $[0,\infty)$ unless $\gamma=0$. It {\em is} possible
for the boundary conditions to be satisfied if $\gamma>0$; for
example, $\P^{m+1}$ is a smooth compactification of the total space of
${\cal O}_{\P^m}(1)$.
\end{remark}

\subsubsection*{Proof of Theorem~\protect\ref{thm:B}}

The proof of Theorem~\ref{thm:B} differs from the proof of
Theorem~\ref{thm:A} in the initial conditions, and consequently in the
choice of momentum interval; $\varphi'(0)=0$ rather than
$\varphi'(0)=2$, so $I=(0,\infty)$ rather than $[0,\infty)$. The
details are almost exactly parallel: The profile of
equation~(\ref{eqn:varphi}) is replaced by 
\begin{equation}
\label{eqn:profile-b}
\varphi(\mum) = 
{2\over{Q(\mum)}}\int_0^\mum (\mum-x)\Bigl(R(x)-c\Bigr)Q(x)\,dx,
\end{equation}
which is positive for $\mum>0$ if $c \ll 0$. An interval ${\rm
J}^\times$, with supremum equal to $c_0^\times$, is defined as
before. The proof that the metrics are complete and live on
$\Delta^\times(L)$ if $c< c_0^\times$ is entirely analogous to the
proof of Lemma~\ref{lemma:complete}. The borderline profile is not
identically zero if $R(\tau)$ is non-constant, and induces a metric
on~$\lx$. 

Exactly as before, the metric is Einstein iff
$$\rho_M=\lambda\omega_M\quad{\rm and}\quad 
{1\over {2Q}}(\varphi Q)'(\mum)=-\lambda\mum,\qquad {\rm so}\quad
\varphi(\mum)=-{{2\lambda}\over{Q(\mum)}}\int_0^\mum xQ(x)\,dx.$$
This function~$\varphi$ is positive on~$(0,\infty)$ iff $\lambda<0$. 
Conversely, if $\rho_M=\lambda\omega_M$ and $c=\lambda(m+1)$ for some
$\lambda<0$, then an argument analogous to the proof of
Lemma~\ref{lemma:corresp} shows that 
$$\Bigl(R(\mum)-c\Bigr)Q(\mum)
=2\lambda{d\over{d\mum}}\Bigl(-\mum Q(\mum)\Bigr),$$
so the associated metric is Einstein. \QED
\medskip

One significant difference between Theorems~\ref{thm:A}
and~\ref{thm:B} appears when $R$~is constant.  In the event $g_M$~is
Ricci-flat (so $R\equiv0$) there is no profile satisfying the
boundary conditions $\varphi(0)=\varphi'(0)=0$ that induces a
scalar-flat metric. In fact, the equation for scalar-flatness is
$(\varphi Q)''=0$, while the initial conditions imply $(\varphi
Q)(0)=(\varphi Q)'(0)=0$, so $\varphi\equiv0$. Even by relaxing the
initial conditions, the only scalar-flat metrics that arise are
uninteresting. To wit, the equation $(\varphi Q)''=0$ leads to the
candidate profile $\varphi(\mum)=(a_0+a_1\mum)/Q(\mum)$. This function
induces a complete metric only if it is everywhere positive (and
perhaps has a removable discontinuity), but in this event the momentum
interval is all of~$\R$. The requirement that $\omega_M-\mum\gamma$ be
a K\"ahler form for all $\mum\in\R$ forces $\gamma=0$, so $Q\equiv1$
and the profile reduces to a positive constant. The induced metric is
a local product of $g_M$ and a flat cylinder of radius equal to the
value of the profile.

\subsubsection*{Bounds on $c_0$}

To facilitate the proof of Lemma~\ref{lemma:finite}, a rational
function~$C$ was introduced whose infimum over $I$ is equal to
$c_0$. Since $C$ has an explicit expression in terms of the curvature
of the horizontal data, $c_0$ can in principle be estimated in terms
of the curvature.  A crude estimate comes from the following simple
observations.  If $R(\mum)-c\geq0$ for all $\mum\geq0$, then the
profile~(\ref{eqn:varphi}) is positive for $\mum>0$, while if
$R(\infty)-c<0$, then the profile is {\em not} always positive. Thus
\begin{equation}
\label{eqn:R-bounds}
\inf_{\mum \geq 0} R(\mum) \leq c_0 \leq R(\infty).
\end{equation}
To estimate the lower bound, pick a point $z$ of~$M$ and
choose an orthonormal basis of $T_z^{1,0}M$ relative to which $\Beta$
is diagonal. Denote by $\beta_\nu$ the eigenvalues of $\Beta$, indexed
so that $\beta_1\leq\cdots\leq\beta_j<0=\beta_{j+1}=\cdots \beta_m$,
and let $\varrho_\nu$ be the diagonal elements of $\varrho$ in the
basis. Then
$$
R(\mum) = \frac{\varrho_1}{1 -\beta_1\mum} + 
\cdots + \frac{\varrho_j}{1 -\beta_j\mum} + \varrho_{j+1} + \cdots
\varrho_m.
$$
Since the coefficients of $R(\mum)$ are by hypothesis independent
of~$z$, certain combinations of the $\varrho_\nu$ are independent
of~$z$.  Indeed suppose that among the $\beta_\nu$ the {\em distinct}
real numbers occuring are $b_1<\cdots< b_{\ell-1} < b_\ell=0$ with
multiplicities $k_1,\ldots, k_\ell$, and let
$$r_1 = \varrho_1+ \cdots+\varrho_{k_1},\quad
r_2 = \varrho_{k_1+1}+ \cdots+\varrho_{k_1+ k_2},\ \ldots,\quad
r_\ell = \varrho_{k_\ell+1}+\cdots+\varrho_m.$$
Then $r_\ell = R(\infty)$ and
$$R(\mum) = \sum_{i=1}^\ell \frac{r_i}{1-b_i \mum} = R(\infty) + 
   \sum_{i=1}^{\ell-1}\frac{r_i}{1-b_i \mum}.$$
Considering each fraction separately, $\displaystyle R(\mum)
\geq R(\infty) + \sum_{i=1}^\ell \min(0,r_i)$, so that
$$R(\infty) + \sum_{i=1}^\ell \min(0,r_i) \leq c_0 \leq R(\infty).$$
Although crude, this is sometimes sufficient to determine $c_0$
exactly. For example, if $\rho_M$ is positive semi-definite, then all
$r_i\geq 0$ and $c_0 = r_\ell=R(\infty)\geq0$.  If in addition
$\gamma$ is negative-definite, then $c_0=0$. In either case, the total
space of~$L$ admits a complete, scalar-flat K\"ahler metric.

In the setting of Theorem~\ref{thm:B}, there is the additional
inequality $c_0^\times \leq R(0)$, necessitated by positivity of the
profile (\ref{eqn:profile-b}) near $\mum=0$. If the base metric has
non-positive Ricci tensor, then $R$ is monotone increasing and
$c_0^\times = R(0) = \displaystyle\sum_{\nu} \varrho_\nu$.

\subsection{Metrics on Vector Bundles}
\label{section:vbundle}

Let $p:(E,h)\to(D,g_D)$ be an Hermitian holomorphic vector bundle of
rank~$n>1$ over a K\"ahler manifold of dimension~$d$. There is a
smooth, globally defined norm squared function~$r$, just as for line
bundles, and the Calabi ansatz has an obvious formulation in this
situation, namely to consider closed $(1,1)$-forms
\begin{equation}
\label{eqn:ansatz2}
\omega=p^*\omega_D+\ddbar F(r)=p^*\omega_D+2\ddbar f(t).
\end{equation}
The machinery of Theorem~\ref{thm:CKS} is easily modified to treat
this case. Complex-analytically, the idea is to blow up the zero
section of~$E$ and pull the form~(\ref{eqn:ansatz2}) back to the total
space of the blow-up. The original horizontal data pull back to
``partially degenerate'' horizontal data since $\omega_D$ has rank~$d$
as an Hermitian form while the exceptional divisor has
dimension~$d+n-1$. To avoid problems caused by this degeneracy, the
form $\omega_D-\gamma$ on~$\P(\Ex)$ is used as a background
metric. The following discussion elaborates on the technicalities
necessary to make this idea work.

\subsubsection*{The tautological bundle, and blowing up}

The complement of the zero section of $E$ is denoted by $\Ex$, and is
regarded either as a complex manifold or as a ``punctured
$\C^n$-bundle'' over~$D$.  There is a free $\C^\times$-action on~$\Ex$
by scalar multiplication, whose quotient is the {\em projectivization}
of~$E$, denoted~$\P(E)$. The complex manifold $M=\P(E)$ is the
total space of a $\P^{n-1}$-bundle $\pi:M\to D$.  The pullback of~$E$
to~$M$ is the fibre product
\begin{equation}
\label{eqn:pullback}
\pi^*E=\{(\zeta,v)\in M\times E\mid \pi(\zeta)=p(v)\}
=M\times_D E
\end{equation}
endowed with the projection onto the first factor, and the {\em
tautological bundle} $\tau_E$ is the line subbundle of $\pi^*E\to M$
whose fibre at a point $\zeta=(z,[v])\in\P(E)=M$ is the line through
$v\in E_z$. Two important observations are:
\begin{itemize}
\item The fibre of~$\pi$ over $z\in D$ is the $(n-1)$-dimensional
projective space $\P(E_z^\times)$, and the restriction of $\tau_E$ to
a fibre $\P(E_z^\times)$ is ${\cal O}(-1)$. The total space
of the latter may be regarded as the blow-up of $\C^n$ at the origin.
\item Projection to the second factor in~(\ref{eqn:pullback})
induces a biholomorphism $\tau_E^\times\simeq\Ex$, and coincides
with the map $\pi:M\to D$ along the zero section of~$\tau_E$.
\end{itemize}
The total space of~$\tau_E$ is obtained from the total space of~$E$ by
{\em blowing up the zero section}. This is a direct generalization of
blowing up a point, and indeed may be regarded as a family of blow-ups
of~$\C^n$, parametrized by points of~$D$.

By abuse of notation, the projection $\tau_E\to M$ is denoted by~$p$.
These spaces and bundles are organized as follows, with superscripts
denoting ranks of bundles or dimensions of spaces:
$$\begin{array}{ccccc}
\pi^*E & \supset & \tau_E  & {\longrightarrow}   &     E^n        \\
\Big\downarrow & & \phantom{p}\Big\downarrow p & &
		\phantom{p}\Big\downarrow p        \\ 
\P(E) & = & M^m          &\stackrel{\pi}{\longrightarrow} & D^d
\end{array}$$
Of course, $n+d=m+1$ since $\Ex\simeq\tau_E^\times$.  Via this
identification, $L=\tau_E$ acquires an Hermitian structure, also
denoted~$h$, from the Hermitian structure of~$E$, and as before
$\gamma=\gamma_1(L,h)$ denotes the curvature form. The norm squared
function $r:\Ex\to(0,\infty)$ may be regarded as a function on~$L$,
and $p^*\gamma=-\ddbar\log r$ as closed $(1,1)$-forms on the total
space of~$L$. Calabi~\cite{Cal0} calls $\gamma$ the ``bi-Hermitian
curvature form'' of~$(E,h)$. The preceeding discussion makes it clear
geometrically that the auxiliary bundle~$L$ arises naturally in the
Calabi ansatz on vector bundles. In spite of the fact that the complex
manifolds $\lx$~and $\Ex$ are biholomorphic, it is best to regard them
as differential-geometrically distinct; the bundle~$\Ex$ is completed
by adding a copy of~$D$ (the level set $\{r=0\}\subset E$), while
$\lx$ is completed by adding a copy of~$M$ (the level set
$\{r=0\}\subset L$). Lemma~\ref{lemma:extend2} below illustrates the
importance of this distinction.

\subsubsection*{Degeneracy in the Calabi ansatz}

Let $p:(E,h)\to(D,g_D)$ be an Hermitian holomorphic vector bundle of
rank~$n$ over a K\"ahler manifold. The closed $(1,1)$-form
$\pi^*\omega_D$ may be denoted by $\omega_D$ for brevity, as in
equation~(\ref{eqn:form2}) below.  Associated to a profile~$\varphi$
are functions $\mom$ and $f$ as in Section~\ref{section:KS}.
Upon blowing up the zero section of~$E$, the closed $(1,1)$-form
$\omega=p^*\omega_D+2\ddbar f(t)$ is immediately written, using
results of Section~\ref{section:KS}, as
\begin{equation}
\label{eqn:form2}
\omega = \varphi(\mum)(dt\wedge d^ct)+p^*\omega_D-\mum p^*\gamma.
\end{equation}

The form $\pi^*\omega_D$ on~$M$ is degenerate along the fibres of
$\pi:M\to D$, while under the hypotheses of Theorem~\ref{thm:C} the
form $\pi^*\omega_D-\gamma$ is a K\"ahler form on~$M$, and (by
Lemma~\ref{lemma:indep} below) the horizontal data
$p:(\tau_E,h)\to(M,\omega_M)$ are $\sigma$-constant.  The vector
bundle data $p:(E,h)\to(D,\omega_D)$ are said to be $\sigma$-constant
in this situation. 
%(There is an obvious definition of
%``$\rho$-constancy'' for vector bundle data.) 
Before defining the
analogues of the functions $Q$~and $R$ for vector bundle data, it is
reassuring to verify that the choice of background metric in the
family $\omega_M(\mum)$ is immaterial.

\begin{lemma}
\label{lemma:indep}
If $\omega_M(\mum)$ is a K\"ahler form for every~$\mum>0$, and if
$p:(L,h)\to\bigl(M,\omega_M(\mum_0)\bigr)$ is $\sigma$-constant for
some $\mum_0>0$, then $p:(L,h)\to\bigl(M,\omega_M(\mum)\bigr)$ is
$\sigma$-constant for every $\mum>0$. 
\end{lemma}

\begin{proof}
Fix $\mum_0>0$, and regard $\omega_M(\mum_0)$ as a background metric,
so that all index raising and lowering is done with respect to this
metric.  Let $\Beta$ denote the curvature endomorphism of $(L,h)$; by
assumption, the eigenvalues $\beta_\nu$ are constant on~$M$, and the
trace of the Ricci form $\rho_M(\mum_0)$ with respect to the K\"ahler
form $\omega_M(\mum)$ is constant for all $\mum>0$. Consider the
family $A(\mum)$ of endomorphisms associated to the closed
$(1,1)$-forms $\omega_M(\mum)$. Then $A(\mum_0)$ is the identity, and
$A'(\mum)=-\Beta$, so
$$A(\mum)={\rm I}-(\mum-\mum_0)\Beta\qquad \hbox{for $\mum\geq0$.}$$
The eigenvalues of~$\Beta$ are $\beta_1,\ldots,\beta_d$, and
$-1/\mum_0$ (the latter of multiplicity~$n-1$).  Thus $A(\mum)$ has
constant eigenvalues for each $\mum\geq0$, $A(0)={\rm I}+\mum_0\Beta$
has rank~$d$, and the endomorphisms $A(\mum)$ have the same
eigenbundles. Set
$$Q(\mum,\mum_0) = \det A(\mum)= {{\mum^{n-1}}\over{\mum_0^{n-1}}}
\prod_{\nu=1}^d\Bigl(1-(\mum-\mum_0)\beta_\nu\Bigr),$$
so that for each $\mum_0>0$, $Q$~is a constant-coefficient polynomial
with a zero of order~$(n-1)$ at $\mum=0$. Then $Q(\mum,\mum_0)$ is the
ratio of the volume forms of $\omega_M(\mum)$ and $\omega_M(\mum_0)$,
and is constant on~$M$ for each fixed $\mum>0$, so the Ricci
forms $\rho_M(\mum)$ and $\rho_M(\mum_0)$ coincide. Consequently, if
the trace of the Ricci form $\rho_M(\mum_0)$ with respect to
$\omega_M(\mum)$ is constant for all $\mum>0$, then the same is true
of the Ricci form~$\rho_M(\mum_1)$ for every $\mum_1>0$.
\end{proof}

Suppose $p:(E,h)\to(D,\omega_D)$ are $\sigma$-constant data, and set
$\omega_M(\mum)=\pi^*\omega_D-\mum\gamma$. The background metric
on~$M$ is taken to be $\omega_M:=\omega_M(1)$, at variance with the
notation used for line bundles. As noted in the proof of
Lemma~\ref{lemma:indep}, the Ricci forms~$\rho_M(\mum)$ do not depend
on~$\mum$, and are denoted simply by~$\rho_M$. By contrast, the Ricci
endomorphisms, defined by
$\varrho_M(\mum):=\omega_M(\mum)^{-1}\rho_M$, {\em do} depend
on~$\mum$, though by hypothesis they all have constant trace. More is
true: For each $\mum>0$, the vertical tangent bundle of~$M$ (namely,
$\ker\pi_*$) is an
eigenbundle of $\varrho_M(\mum)$ with eigenvalue~$n$ (of
multiplicity~$n-1$). To see this, observe that the global inner
product $\langle\rho_M,\gamma\rangle$ is constant on~$M$ since
$\langle\rho_M,\omega_M(\mum)\rangle$ is constant on~$M$ for each
$\mum>0$. On a fibre of~$\pi$, which is a projective space, $-\gamma$
is a K\"ahler form, and since the trace of~$\rho_M$ restricted to the
fibre is constant, the restriction must be a multiple of~$-\gamma$,
which implies the Ricci endomorphism is a multiple of the identity on
each fibre. For cohomological reasons, the eigenvalue is~$n$.

The functions $Q$~and $R$ are defined with respect to the background
metric~$\omega_M$ by
\begin{eqnarray}
Q(\mum) & = & \det A(\mum) \quad = \quad \mum^{n-1}\prod_{\nu=1}^d
(1+\beta_\nu-\mum\beta_\nu)
\quad =:\quad \mum^{n-1}\,Q_0(\mum), \nonumber \\
&& \label{eqn:QRdef2} \\
R(\mum) & = & \tr\Bigl(A(\mum)^{-1}\varrho\Bigr)
\quad = \quad {{n(n-1)}\over\mum}+ {\rm smooth}.
\nonumber
\end{eqnarray}
Note that~$Q$ has a zero of order~$(n-1)$ and~$R$ has a simple pole at
$\mum=0$. 

\subsubsection*{Proof of Theorem~\protect\ref{thm:C}}

As in Section~\ref{section:KS}, a profile~$\varphi$ induces a
bundle-adapted metric on~$E$ whose scalar curvature is
\begin{equation}
\label{eqn:scalar2}
\sigma_\varphi(\mum) 
= \left(R-{1\over {2Q}}(\varphi Q)''\right)(\mum).
\end{equation}
The function $\sigma_\varphi$ generally has a pole at
$\mum=0$. Interestingly, the scalar curvature is bounded near the zero
section if, {\em and only if}, the metric extends over the zero
section: 

\begin{lemma}
\label{lemma:extend2}
If the metric~{\rm(\ref{eqn:ansatz2})} is complete, then the function
$\sigma_\varphi$ in equation~{\rm(\ref{eqn:scalar2})} has a removable
singularity at $\mum=0$ if and only if the profile satisfies the
boundary conditions $$\varphi(0)=0,\qquad\varphi'(0)=2,$$ if and only
if the metric extends over the zero section of $E$.
\end{lemma}

\begin{proof}
Completeness of the metric dictates that $\varphi(0)=0$. Writing
$\varphi(\mum)=\varphi'(0)\,\mum + O(\mum^2)$ near $\mum=0$---so that
$(\varphi Q)(\mum)=\varphi'(0)\mum^n\, Q_0(\mum)+O(\mum^{n+1})$---and
differentiating twice gives
$$\sigma_\varphi(\mum)
={{n(n-1)}\over{\mum}}\left(1-{{\varphi'(0)}\over2}\right) + O(1),$$
so the singularity at $\mum=0$ is removable if and only if
$\varphi'(0)=2$.  The remaining assertions are proven exactly as in
the line bundle case.
\end{proof}

Lemma~\ref{lemma:extend2} says the Calabi ansatz~(\ref{eqn:ansatz2})
does not give rise to a complete K\"ahler metric of constant scalar
curvature on the punctured disk subbundle of~$E$. Of course, the
punctured disk subbundle of~$E$ is biholomorphic to the punctured disk
subbundle of the line bundle~$\tau_E$; the difference between
Theorem~\ref{thm:B} and the present situation is that here the forms
$\omega_M(\mum)$ drop in rank at $\mum=0$.

The differential equation obtained by setting (\ref{eqn:scalar2})
equal to $c$ has a regular singular point at $\mum=0$, and the general
solution is
$$\varphi(\mum) = {2\over{Q(\mum)}}\left(\alpha_0+\alpha_1\,\mum 
+ \int_0^\mum\Bigl(R(x)-c\Bigr)(\mum-x)Q(x)\,dx\right).$$
Because $Q$ has a zero of order~$n-1$ at $\mum=0$, $\varphi(0)=0$
forces $\alpha_0=\alpha_1=0$, so the purported profile for a
bundle-adapted metric of scalar curvature~$c$ on the disk subbundle of
$E$ is 
\begin{equation}
\label{eqn:varphi2}
\varphi(\mum) = 
{2\over{Q(\mum)}}\int_0^\mum\Bigl(R(x)-c\Bigr)(\mum-x)Q(x)\,dx.
\end{equation}
This profile satisfies the boundary conditions $\varphi(0)=0$ and
$\varphi'(0)=2$. The first results from two applications of
l'H\^{o}pital's rule. To obtain the latter, differentiate $\varphi Q$
and solve for $\varphi'$ to get
$$\varphi'(\mum)=-\varphi(\mum){{Q'}\over{Q}}(\mum)
+ {2\over Q}\int_0^\mum\Bigl(R(x)-c\Bigr)Q(x)\,dx.$$
Writing $\varphi(\mum)=\varphi'(0)\mum+O(\mum^2)$ and using
l'H\^{o}pital's rule gives 
$$\varphi'(0)= \lim_{\mum\to 0}\left(-\varphi'(0)(n-1) + O(\mum) +
2{{\Bigl(R(x)-c\Bigr)Q(\mum)}\over {Q'(\mum)}}\right) =
-\varphi'(0)(n-1) + 2n,$$ 
so $\varphi'(0)=2$ as claimed. Most of the remaining points are
checked exactly as in the proof of Theorem~\ref{thm:A}. Specifically,
the set ${\rm J}=\{c\mid \varphi(\mum)>0\mbox{ for all $\mum>0$} \}$
is a non-empty interval, bounded above, and $c_0:=\sup \rm J$. The
profile~(\ref{eqn:varphi2}) is positive for all $\mum>0$ if $c<c_0$,
and has at most quadratic growth, while if $c=c_0$ the induced metric
is complete (though possibly of finite fibre volume).  When $c>c_0$,
there is a first positive root~$b$, and three possibilities occur:
\begin{itemize}
\item $\varphi'(b)=0$, so the metric is complete and of finite volume;
\item $\varphi'(b)=-2$, so the metric extends smoothly to
  $\P(E\oplus{\cal O}_D)$, which is the $\P^n$-bundle over~$D$
  obtained from~$E$ by adding a divisor at infinity;
\item $\varphi'(b)\neq 0,-2$, so the metric is incomplete and has no
  smooth extension to a larger manifold.
\end{itemize}
As in Lemma~\ref{lemma:finite}, there are at most finitely many values
of $c>c_0$ such that one of the first two possibilities holds.

The analysis of the Einstein condition is very similar in this case to
the previous two (see Lemma~\ref{lemma:corresp}). One finds that
$\rho_\varphi=\lambda\omega_\varphi$ if and only if
$$\lambda\pi^*\omega_D=\rho_M+n\,\gamma\qquad{\rm and}\qquad
{1\over{2Q}}(\varphi Q)'(\mum)=n-\lambda\mum.$$
The second equation is a differential equation which leads to the
profile whose associated metric has scalar curvature
$\lambda(n+d)$. As a consistency check, observe that the form
$\rho_M+n\,\gamma$ is cohomologically degenerate along the fibres
of~$\pi$, since $\rho_M$ restricts to the curvature of a Fubini-Study
metric, while $\gamma$ restricts to the negative generator of second
cohomology. 

\subsubsection*{Collapsing metrics along the zero section}

Theorem~\ref{thm:C} indicates that the momentum
construction provides a differential-geometric framework in which to
view (partial) collapsing of divisors, provided the normal bundle
satisfies curvature hypotheses.  The following remarks explain this
idea in more detail. Let $p:(E,h)\to(D,g_D)$ be an Hermitian vector
bundle satisfying the hypotheses of Theorem~\ref{thm:C}. For each
$a>0$, the horizontal data $p:(L,h)\to M$ equipped with the base
metric $\omega_M(a)=\pi^*\omega_D-a\gamma$ are $\sigma$-constant, and
the ``critical'' scalar curvature constants $c_0(a)$ are uniformly
bounded below for~$a$ in some interval to the right of~$0$. If
$c<c_0(a)$ for all such~$a$, then the momentum construction yields a
family of metrics on the total space of~$L=\tau_E$, each complete and
of scalar curvature~$c$. As $a\to0$, these metrics converge to a
metric on the total space of~$E$, i.e.\ to a closed $(1,1)$-form
on~$L$ whose restriction to the zero section has rank~$d$ and is
otherwise positive-definite. The fibres of $M=\P(E)\to D$ shrink
homothetically while the curvature of the fibre metric
approaches~$-\infty$ along the zero section in such a way that the
scalar curvature of the total space stays constant.

\section{Limitations, Examples, and Literature Survey}
\label{section:examples}
\setcounter{equation}{0}
\setcounter{figure}{0}

In this section we shall address the following questions:
\begin{enumerate}
\item What is the role of the assumption $I =(0,\infty)$ or
$[0,\infty)$, made throughout Section~\ref{section:inf-vol}, and what
do the methods say about metrics of finite (fibre) volume,
particularly compact metrics?

\item To what extent can the assumption of $\sigma$-constancy be
relaxed? Are there complete metrics of constant scalar curvature that
may be constructed via the Calabi ansatz, but which are not
encompassed in the existence theorems of
Section~\ref{section:inf-vol}?

\item How explicitly are the examples described? Are there any
``genuinely'' new examples?
\end{enumerate}
We shall also give a brief survey of related literature in
Section~4.5. 

Question~1 is treated in Sections 4.1~and 4.4.  As stated, Question~2
is somewhat heuristic, reflecting our lack of a completely definitive
answer. Partial results, including the proof of
Theorem~\ref{thm:constant}, are collected in Section~\ref{scope}.
Examples of $\sigma$-constant data are given in
Section~\ref{sec:examples}.  These provide an answer to Question~3; it
emerges that there is greater freedom in choosing the K\"ahler class
of the base metric~$g_M$ than there is in prior works. On the question
of explicitness, the formulae of Section~\ref{section:inf-vol} show
that $\varphi$ is a rational function whose coefficients depend
explicitly upon the Ricci endomorphism~$\varrho$ and the curvature
endomorphism~$\Beta$; see also Section~\ref{sec:examples}
below. Section~\ref{section:KS} exhibits some important geometric
invariants (Ricci/scalar curvature, and Laplace operator) of the
metrics very simply and explicitly in momentum coordinates.  The price
to be paid is that the transformation back to holomorphic coordinates
involves inversion of an integral.

\subsection{Scope of the Momentum Construction}
\label{scope}

This subsection describes our attempt to understand the limitations of
the Calabi ansatz. The general construction begins with horizontal
data---an Hermitian holomorphic line bundle over a complete K\"ahler
manifold, and an interval of real numbers---and associates to each
profile~$\varphi$ a K\"ahler metric~$g_\varphi$. There are two
intervals in the Calabi ansatz: The momentum interval~$I$, which is
closely related to the area of the fibre metric; and the $r$~interval
defining the invariant subbundle on which the metric is defined. That
these intervals are essentially independent is clear from
Table~\ref{table-big}.

\subsubsection*{Normalization of the momentum interval}

Begin with horizontal data $\{p:(L,h)\to(M,g_M),\ I\}$. The horizontal
forms $\omega_M-\mum\,\gamma$ are assumed to be positive and complete
for all $\mum\in I$, so if $I=\R$, then $\gamma=0$, and every profile
gives rise to a local product metric. In all other cases, $I$ may be
normalized without loss of generality; in words, every bundle-adapted
metric is isometric to a bundle-adapted metric having momentum
interval~$I$ equal to $\R$, or else having $\inf I=0$:

\begin{lemma}
\label{lemma:interval}
If $\{p':(L',h')\to(M,g_M'),\ I'\}$ are horizontal data with
$I'\neq\R$ and if $\psi$ is a profile on~$I'$ inducing a
metric~$g_\psi$, then there exist horizontal data
$\{p:(L,h)\to(M,g_M),\ I\}$, with $\inf I=0$, and a profile~$\varphi$
on~$I$ such that $g_\varphi$ and~$g_\psi$ are isometric.
\end{lemma}

\begin{proof}
Suppose first that $I'$ is bounded below, and set $a=\inf I'$; thus
$\omega_M'-a\gamma'$ is a K\"ahler form on~$M$. The `translated' data 
$$(L,h)=(L',h'),\qquad \omega_M=\omega_M'-a\gamma',\qquad
I=I'-a,\qquad \varphi(\mum)=\psi(\mum-a),$$
have the advertized properties.

Now suppose $I'\neq\R$ is not bounded below. By a translation argument
analogous to that just given, it may be assumed that $\sup I'=0$. 
Consider the `inverted' data
$$(L,h)=({L'}^*,h^{-1}),\qquad \omega_M=\omega_M',\qquad
I=-I',\qquad \varphi(\mum)=\psi(-\mum).$$
The $t$~integral acquires a sign change, which corresponds to
the inversion map $\iota:(L,h) \to (L^*,h^{-1})$, given
locally by $z^0 \mapsto 1/z^0$ in a line bundle chart. It is
straightforward to verify $g_\varphi=\iota^*g_\psi$.  
\end{proof}

\subsubsection*{Level sets of $\sigma_\varphi$}

According to Theorem~\ref{thm:CKS} and equation~(\ref{sigma2}), the
scalar curvature~$\sigma_\varphi$ is given in terms of horizontal data
and the profile by 
\begin{equation}
\label{eqn:scal}
\sigma_\varphi = \sigma_M(\mum)-{1\over{2Q}}
{{\del^2}\over{\del\mum^2}}(Q\varphi)(\mum),
\qquad \sigma_M(\mum)=R(\mum)-\lapl_{\omega_M(\mum)}\log Q;
\end{equation}
Our aim has been to find data (bundle, metrics, interval, and profile)
such that $\sigma_\varphi$ is constant. As noted in the introduction,
it cannot be expected that the level sets of~$\sigma_\varphi$---which
depend on the curvature of $(L,h)$---will coincide with the level sets
of~$\mum$ for an arbitrary profile. Similarly, it is not to be
expected that the Calabi ansatz yields a metric of constant scalar
curvature for arbitrary horizontal data. The goal is thus to determine
curvature conditions on the horizontal data that are necessary and
sufficient for existence of a profile inducing a metric of constant
scalar curvature (or perhaps satisfying some weaker condition, such as
having scalar curvature depending only upon~$\mum$). The condition of
$\sigma$-constancy is sufficient, by the existence theorems of the
preceeding section. The aim of this subsection (unfortunately not
completely realized) is to determine the extent to which
$\sigma$-constancy is a necessary condition.

Equation~(\ref{eqn:scal}) holds for all profiles and horizontal
data. Regard specification of the scalar curvature~(\ref{eqn:scal}) as
a family of ODEs on~$I$, parametrized by points of~$M$. This yields a
family of solutions, viewed as a function $\varphi:I\times M\to\R$,
and the question is to determine when this function depends only
on~$\mum\in I$. 
For the moment, assume there is a (smooth) function~$\sigma$ on~$I$ so
that $\sigma_\varphi=\sigma(\mum)$. Then 
\begin{equation}
\label{eqn:profile-family}
\varphi(\mum,z)={1\over{Q(\mum,z)}} 
\left[{{\del\varphi}\over{\del\mum}}(0,z)\,\mum +
2\int_0^\mum(\mum-x)\Bigl(\sigma_M(x,z)-\sigma(x)\Bigr)
Q(x,z)\,dx\right].
\end{equation}
Our general goal is to make minimal hypotheses on~$\sigma$ (e.g.\
$\sigma$~is real-analytic) and then to prove that if the right-hand
side of~(\ref{eqn:profile-family}) depends only on~$\mum$ for a
particular choice of~$\sigma$, then the horizontal data are
$\sigma$-constant. 

\subsubsection*{Proof of Theorem~\protect{\ref{thm:constant}}, and
related remarks}

Fix horizontal data $\{p:(L,h)\to(M,g_M),\ I\}$. Suppose there exist
real-analytic functions $\sigma_1$~and $\sigma_2$ on~$I$ so that the
corresponding profiles $\varphi_1$~and $\varphi_2$ given
by~(\ref{eqn:profile-family}) depend only on~$\mum$ and agree at
$\mum=0$. Then the function 
$$(\varphi_2-\varphi_1)(\mum)={1\over{Q(\mum,z)}}\left[
\Bigl(\varphi_2'(0)-\varphi_1'(0)\Bigr)\mum + \int_0^\mum(\mum-x)
\Bigl(\sigma_1(x)-\sigma_2(x)\Bigr)Q(x,z)\,dx\right]$$
depends only on~$\mum$ and vanishes at~$0$. Let $N(\mum,z)$ denote the
term in square 
brackets, and let $\psi(\mum)=N(\mum,z)/Q(\mum,z)$. The aim is to show
that the horizontal data are $\sigma$-constant, namely that
$Q(\mum,z)$ depends only on~$\mum$ and the horizontal metric
$g_M(\mum)$ has constant scalar curvature for each $\mum\in I$.
Because $\sigma_\varphi$ is {\em assumed}\/ to depend only
upon~$\mum$ (for $\varphi=\varphi_1$, say), it suffices to show
that~$Q$ depends only on~$\mum$, since
$\sigma_M(\mum)=\sigma_\varphi+(1/2Q)(\varphi Q)''(\mum)$. To this
end, it is enough to show that $N$ depends only on~$\mum$ since
$$N''(\mum,z)=\Bigl(\sigma_1(\mum)-\sigma_2(\mum)\Bigr)Q(\mum,z).$$ 

The function~$\psi$ is real-analytic and not identically zero (unless
$\varphi_1\equiv\varphi_2$). Put $\alpha=\varphi_2'(0)-\varphi_1'(0)$, 
write $\sigma_1(\mum)-\sigma_2(\mum) = S(\mum)$, and consider the
initial-value problem  
\begin{equation}
\label{eqn:IVP}
y''(\mum)-{{S(\mum)}\over{\psi(\mum)}}y(\mum)=0,
\qquad y(0)=0,\quad y'(0)=\alpha.
\end{equation}
The idea is that for each $z\in M$, the function $y=N(\ ,z)$
satisfies~(\ref{eqn:IVP}), so by uniqueness it follows that $N$ is
independent of~$z$. The ODE has a singular point at $\mum=0$, however,
so the desired conclusion requires a modicum of additional work.

If $\sigma_1(\mum)\equiv \sigma_2(\mum)$, then it is clear that
$Q(\mum,z)$ does not depend on~$z$, and the proof is
finished. Otherwise, $S(\mum)=S_0\mum^\ell+O(\mum^{\ell+1})$ with
$S_0\neq0$. Since $Q(0,z)=1$ for all~$z$, $Q(\mum,z)S(\mum)=S_0
\mum^\ell +O(\mum^{\ell+1})$, and integrating twice gives
$$N(\mum,z)=\alpha\mum + {{S_0}\over{(\ell+1)(\ell+2)}}\mum^{\ell+2}
+O(\mum^{\ell+3})$$
near $\mum=0$. According to whether $\alpha\neq0$ or $\alpha=0$,
the coefficient $S/\psi$ in~(\ref{eqn:IVP}) is given by
$$\mum^{\ell-1}\Bigl({{S_0}\over\alpha}+O(\mum)\Bigr)\quad{\rm or}
\quad {{(\ell+1)(\ell+2)}\over{\mum^2}}\Bigl(1+O(\mum)\Bigr),$$
cf.~(\ref{eqn:coeff}) below, so the ODE has at worst a regular
singular point at $\mum=0$.  Fix a point of $z\in M$, and expand $N$
and $S/\psi$ as power series about $\mum=0$: 
$$N(\mum,z)=\sum_{n=0}^\infty a_n(z)\mum^n,\qquad\qquad
{{S(\mum)}\over{\psi(\mum)}} 
= {1\over{\mum^2}}\sum_{k=0}^\infty b_k\mum^k.$$
The initial conditions in~(\ref{eqn:IVP}) imply $a_0=0$ and
$a_1=\alpha$. Setting $y=N(\ ,z)$ in the ODE and equating coefficients
gives, suppressing the dependence of~$a_n$ on~$z$, 
\begin{equation}
\label{eqn:coeff}
a_1 b_0 = 0,\qquad
\Bigl[(n+2)(n+1)-b_0\Bigr]a_{n+2}=\sum_{k=0}^{n+1}a_k b_{n-k}
\quad{\rm for}\quad n\geq0.
\end{equation}
It is convenient to treat the cases $\alpha\neq0$ and $\alpha=0$
separately. 
\smallskip

\noindent ($\alpha\neq0$) \quad The function $S/\psi$ has a simple
pole at $\mum=0$, so $b_0=0$. By induction on~$n$, the recursion
relation~(\ref{eqn:coeff}) expresses the coefficients $a_n(z)$ as
linear combinations of $a_0$, $a_1$, and the $b_k$. Since the latter
are assumed constant, the coefficients $a_n(z)$ do not depend on~$z$,
proving that $N(\mum,z)$ depends only on~$\mum$.
\smallskip

\noindent ($\alpha=0$) \quad In this case, $N$ vanishes to order
$(\ell+2)$ at $\mum=0$, so $a_n(z)\equiv0$ for $0\leq n\leq\ell+1$,
while $a_{\ell+2}(z)\equiv S_0/(\ell+1)(\ell+2)$ for all $z\in M$;
thus $b_0=(\ell+2)(\ell+1)$. As before, the
recursion~(\ref{eqn:coeff}) expresses the coefficients $a_n(z)$,
$n>\ell+2$, in terms of the constants $a_{\ell+2}$ and $b_k$, so 
$N(\mum,z)$ depends only on~$\mum$.

This completes the proof Theorem~\ref{thm:constant}. \QED
\medskip

Even the case $\sigma_1=\sigma_2=c$ is of geometric interest: If
$\{p,\ I\}$ are 
horizontal data, and there exist fibre-complete, bundle-adapted
metrics of scalar curvature~$c$ in $\Delta(L)$ (so that
$\varphi_1'(0)=2$) and $\Delta^\times(L)$ (i.e.\ $\varphi_2'(0)=0$),
then the data are $\sigma$-constant, cf.~Theorems \ref{thm:A} and
\ref{thm:B}. It is worth emphasizing that Theorem~\ref{thm:constant}
is local; nothing is assumed about completeness of the base or fibre
metric. Though the result is geometrically significant, the proof is
elementary because of Theorem~\ref{thm:CKS} and the assumption of
existence of {\em two}\/ profiles whose scalar curvature is
independent of~$\mum$; taking the difference $\varphi_2-\varphi_1$
cancels the relatively complicated (but profile-independent) term
$\sigma_M(\mum)$. 

By Theorem~\ref{thm:constant}, it is reasonable to assert that the
Calabi ansatz is exhausted, so far as families of K\"ahler metrics of
constant scalar curvature are concerned. More precisely, if $\{p,\
I\}$ is {\em not}\/ $\sigma$-constant, then there is at most one
metric of constant scalar curvature that arises from the Calabi ansatz
in the disk bundle~$\Delta(L)$.  For the remainder of this section, we
consider the more difficult question of deducing curvature properties
from existence of a {\em single}\/ profile.

\begin{problem}  
Let $\{p,\ I\}$ be horizontal data, and suppose there exists a
profile~$\varphi$ depending only on~$\mum$ such that
$\sigma_\varphi$ is constant. Are the horizontal data necessarily
$\sigma$-constant? 
\end{problem}

\noindent We are not presently able to answer this question, though we
believe the answer is ``yes.'' While falling short of a proof, our
calculations provide evidence and yield some suggestive partial
results. 

The approach is to draw conclusions about the curvature of horizontal
data under the assumption that there exists a germ of a profile at
$\mum=0$, satisfying $\varphi(0)=0$ and $\varphi'(0)=2$, and inducing
a bundle-adapted metric whose scalar curvature depends only on~$\mum$.
This investigation only involves consideration of rational functions
defined on $I\times U$, with $U$ a neighbourhood of the origin
in~$\C^m$. Solutions of ODEs defining momentum profiles are considered
in the `non-geometric regime' where $\mum<0$. Nothing is being
asserted about bundle-adapted K\"ahler metrics for $\mum<0$, of
course; what is being used is real-analyticity of~$\varphi$ in~$\mum$,
and the fact that if such a function $\varphi(\mum,z)$ depends only
on~$\mum$ for $\mum\geq0$, then the same is true for $\mum<0$.

\begin{proposition} 
Fix horizontal data $\{p:(L,h)\to(M,g_M),\ I\}$ and a polynomial
function~$\sigma$.  If there exists a profile~$\varphi(\mum)$ such
that $\sigma_\varphi = \sigma(\mum)$, then $\varphi$ is a {\em
rational} function. 
\end{proposition}

\begin{proof}  
Fix $z\in M$, and let $\varphi(\mum,z)$ be the function defined
by~(\ref{eqn:profile-family}). It will be shown that a
non-rational term in~$\varphi$ has a logarithmic singularity at a root
of~$Q$, and that this root must depend on~$z$, so~$\varphi$ depends
on~$z$ as well. These logarithmic terms potentially arise from 
$\lapl_{\omega_M(\mum)}\log Q$, which enters via
$\sigma_M(\mum)=R(\mum)-\lapl_{\omega_M(\mum)}\log Q$. 

It is enough to work in a coordinate neighbourhood~$U$ around
$z\in M$. Let $\{\beta_\nu(z)\}_{\nu=1}^m$ denote the ({\em
a~priori}\/ non-constant) eigenvalues of the curvature
endomorphism~$\Beta$ in~$U$, and let $b_1<\cdots<b_\ell\leq0$ be the
{\em distinct}\/ eigenvalues, of multiplicity~$k_i$; without loss of
generality, $k_i$ may be assumed constant throughout~$U$.  Then
$$Q(\mum,z)=\prod_{\nu=1}^m\Bigl(1-\mum \beta_\nu(z)\Bigr)
=\prod_{i=1}^\ell\Bigl(1-\mum b_i(z)\Bigr)^{k_i},\qquad
\log Q(\mum,z)=\sum_{i=1}^\ell k_i\log\Bigl(1-\mum b_i(z)\Bigr).$$
Choose a unitary eigenframe $\{\e_\nu\}_{\nu=1}^m$ in~$U$, and let
$(z^\nu)$ be complex local coordinates such that for $1\leq\nu\leq m$,
$z^\nu=0$ and $\del_\nu=\e_\nu$ at~$z$. Differentiating $\log
Q(\mum,z)$, using subscripts (after a comma) to denote partial
derivatives, gives 
\begin{equation}
\label{ddlogq}
\del_\lambda\delbar_\mu \log Q = \sum_{i=1}^\ell\left(
\frac{k_i\mum^2 b_{i,\lambda}b_{i,\bar\mu}}{(1 - \mum b_i)^2}
-\frac{k_i\mum b_{i,\lambda\bar\mu}}{(1-\mum b_i)}\right).
\end{equation}
Let $\{\theta_\nu\}$ be the basis of $\bigwedge_{\bf0}^{(1,0)}$ dual
to $\{\e_\nu\}$, so along the fibre at~$\bf0$ the metric and its
inverse are given by
$$
\omega_M(\mum) 
= \sum_{\nu=1}^m (1-\beta_\nu\mum)\, \theta_\nu\bar{\theta}_\nu,
\qquad
\omega_M(\mum)^{-1} 
= \sum_{\nu=1}^m \frac{e_\nu\bar{e}_\nu}{(1-\beta_\nu\mum)}.
$$ 
Taking the trace of~(\ref{ddlogq}) with respect to
$\omega_M(\mum)^{-1}$, and using the fact that $\{\e_\nu\}$ is a
coordinate basis at~$\bf0$, 
\begin{equation}
\label{eqn:laplogq}
\lapl_{\omega_M(\mum)}\log Q(\mum,z)
= \sum_{i=1}^\ell\sum_{\lambda=1}^m
\left(\frac{k_i\mum^2 |b_{i,\lambda}|^2}
{(1-\mum \beta_\lambda)(1 - \mum b_i)^2}
-\frac{k_i\mum b_{i,\lambda\lambda}}
{(1-\mum \beta_\lambda)(1-\mum b_i)}\right)
\end{equation}
at~$\bf0$. When the right-hand side is expanded in partial fractions,
the coefficients of the principal part at $\mum = 1/b_i$ are linear
combinations of derivatives of~$b_i$.  In particular, if a  principal
part appears, then $b_i$ is non-constant.

The integrand in~(\ref{eqn:profile-family}) is rational for each $z\in
M$, and the only non-polynomial terms come
from~(\ref{eqn:laplogq}). In other words, $\varphi(\mum,z)$ is the
integral of a rational function, and any non-rational (i.e.\
logarithmic) terms arise from a principal part. But as just observed,
a principal part is non-vanishing only when some~$b_i$ is
non-constant, and this implies $\varphi$ depends on~$z$.
\end{proof}

Under similar hypotheses, the calculations in the proof can be used to
read off other geometrically interesting consequences. When the
curvature endomorphism~$\Beta$ is a multiple of the identity,
examination of the highest-order pole in~(\ref{eqn:laplogq}) implies
$\sigma$-constancy:

\begin{proposition} 
Let $\{p,\ I\}$ be horizontal data such that $\Beta(z)=\beta(z){\rm
I}$ for some smooth function~$\beta$. Assume there is a
profile~$\varphi(\mum)$ inducing a metric whose scalar curvature is
$\sigma(\mum)$ for some polynomial~$\sigma$. Then the data are
$\sigma$-constant.
\end{proposition}

Finally, observe that the non-constant curvature eigenvalues are
tightly constrained by the assumptions made above.

\begin{proposition}
Let $\{p,\ I\}$ be horizontal data, and let $\{b_i(z)\}_{i=1}^\ell$
denote the distinct eigenvalues of the curvature endomorphism.
Suppose there exists a profile, depending only on~$\mum$ and
inducing a metric whose scalar curvature depends only on~$\mum$. If
$b_i(z)$ is  non-constant, then for every $z\in M$, $(1-\mum
b_i)^{k_i}$  divides the term in square
brackets in~{\rm(\ref{eqn:profile-family})}. 
\end{proposition}

\begin{proof}
By hypothesis, $\varphi$---and in particular each pole---is
independent of~$z$. So if $(1-\mum b_i)^{k_i}$ does not divide the
numerator (the term in square brackets) in~(\ref{eqn:profile-family})
then $b_i$ is independent of~$z$. 
\end{proof}
The condition that the numerator should have such a factor is a strong
constraint on the way in which $b_i(z)$ can vary.  Unfortunately, it
is a constraint whose consequences seem hard to express in simple
fashion.  

\subsection{Line Bundles}
\label{sec:examples}

This section describes some examples of data satisfying the hypotheses
of Theorems \ref{thm:A}~and \ref{thm:B}. The ``atomic'' examples are
suitable line bundles over Hodge manifolds of constant scalar
curvature; even when the base is a curve, interesting (and seemingly
new) metrics are obtained.

\subsubsection*{Scalar-flat metrics on complex surfaces}

It is instructive to see what emerges from the momentum construction
when the (complex) dimension of the base $M$ is $1$. (The case $\dim M
=0$ was dealt with in Section~\ref{section:complete}, see
Table~\ref{table-big}.) 

The conditions of $\sigma$-constancy are satisfied iff $\rho_M =
\lambda\omega_M$ and $\gamma = \beta\omega_M$, for constants $\lambda$
and $\beta$. Then $Q(\mum) = 1 -\beta\mum$ and $R(\mum) =
\lambda(1-\beta\mum)^{-1}$. The profiles with $\varphi(0)=0$ which
yield a metric of constant scalar curvature~$c$ are given, with their
derivative at~0, by
\begin{equation} \label{prof1}
\varphi(\mum) = \frac{2\mum + (\lambda - c)\mum^2 + c\beta\mum^3/3}
{1- \beta\mum},\qquad \varphi'(0)=2;
\end{equation}
or
\begin{equation} \label{prof2}
\varphi(\mum) = \frac{(\lambda - c)\mum^2 + c\beta\mum^3/3}
{1- \beta\mum},\qquad \varphi'(0)=0.
\end{equation}
It is easy to analyze the choices of $\beta$, $\lambda$, and~$c$ that
give rise to positive profiles and hence to complete metrics. However,
the most interesting case is that where $c=0$, for then the metric is
anti-self-dual in the sense of $4$-dimensional conformal geometry. In
addition to recovering a number of known examples, we find new metrics
on $\C^2$ and its quotients by $\Z$ and $\Z \oplus \Z$.

The first case to consider is when $\lambda>0$, so $M=\P^1$ with a
Fubini-Study (round) metric, normalized so that $\lambda=1$. In this
case, if $L = \O(-k)$, then $\beta = - k/2$; the factor of~2 arises
because the canonical bundle is $\O(-2)$. Equation~(\ref{prof1})
takes the form 
$$
\varphi_k(\mum) = \frac{2\mum + \mum^2}{1+k\mum/2},
$$
which is clearly positive on $(0,\infty)$ if $k \geq0$. The growth is
linear at~$\infty$, so the corresponding metric~$\omega_k$ lives
on the line bundle (not some disk subbundle). 

By Theorem~\ref{thm:A}, $\omega_k$ is Einstein iff $k=2$, in which
case  it is an example of the Ricci-flat metric on~$T^*\P^d$ found by
Calabi~\cite{Cal0}, also known as the Eguchi-Hanson graviton.  If
$k=1$ the metric is the Burns metric on the blow-up of~$\C^2$ at the
origin \cite{Leb2,Lebrun}.  If $k \geq 3$ the metrics are those found
by LeBrun in \cite{Leb2}.

In the case $\lambda <0$, equation~(\ref{prof1}) never yields
complete metrics; the corresponding profiles are all negative
for large~$\mum$, but if $b$ is the first positive zero then
$\varphi'(b) = -2$ iff $\beta=0$.  This case does, however, yield
compact extremal K\"ahler metrics, see~\cite{CTF}.

Finally, if $\lambda=0$---so the universal cover of $M$ is $\C$---then
(\ref{prof1}) reduces to
\begin{equation}
\label{eqn:flat-profile}
\varphi(\mum) = \frac{2\mum}{1 -\beta\mum}.
\end{equation}
This profile yields a complete metric on~$L=\O_\C$, whose total space
is~$\C^2$, provided $\beta < 0$.  These complete, scalar-flat K\"ahler
metrics on~$\C^2$ are not Einstein, and appear to be new. The formula
$h_\beta := \exp(\beta(\Im z)^2)$ defines an Hermitian metric of
constant, negative curvature; indeed, $\gamma(\O_\C,h_\beta) =
\beta\omega$, where $\omega = (\sqrt{-1}/2)dz\wedge d\bar{z}$ is the
standard K\"ahler form.

If the Hermitian metric~$h$ is suitably translation-invariant, then
the profile~(\ref{eqn:flat-profile}) gives rise to complete metrics of
zero scalar curvature on $(\C/\Z)\times\C$ and on the total space of a
line bundle of negative degree over an elliptic curve. If $\beta$ is
integral, $h_\beta$ is an example of such a $\Z\oplus \Z$-invariant
metric.

\subsubsection*{Bundles over product manifolds}
\label{section:4.1}

The basic example of $\sigma$-constant horizontal data comes from a
combination of some well-known results:

\begin{lemma}
\label{lemma:atom}
Let $(M,g)$ be a Hodge manifold  of constant scalar curvature. Then
there exists an Hermitian line bundle $p:(L,h)\to(M,g)$ such that
the data $\{p,[0,\infty)\}$ are $\sigma$-constant.
\end{lemma}

\begin{proof}
By the Lefschetz Theorem on $(1,1)$-classes, there exists a
holomorphic line bundle~$L$ whose first Chern class is $[-\omega]$. To
see there is an Hermitian structure~$h$ with $\gamma(L,h)
=-2\pi\omega$, start with an arbitrary Hermitian structure $h_0$ and
let $\gamma_0$ be the curvature form. By the Hodge Theorem, there
exists a unique smooth, real-valued function~$u$ satisfying
$$\gamma-\gamma_0=\ddbar u,\qquad \int_M u\,\dvol=0.$$ Put
$h=e^{-u}h_0$.  For this Hermitian structure, the horizontal forms
$\omega_M(\mum)$ are positive---indeed, are homothetic
to~$\omega_M$---hence have constant scalar curvature for all
$\mum\geq0$. Finally, the curvature endomorphism is a scalar multiple
of the identity, in particular has constant eigenvalues.
\end{proof}

There are trivial improvements on the statement; the K\"ahler form
need only be {\em homothetic}\/ to an integral form, and every
positive power of~$L$ admits a suitable Hermitian structure.
It is also clear that there are examples of $\sigma$-constant
horizontal data over certain non-compact manifolds, such as Hermitian
symmetric spaces. 

\begin{remark}\label{rem:iterate}
It is tempting to `iterate' the momentum construction, using a
constant scalar curvature metric~$g_\varphi$ on a disk
bundle~$\Delta(L)$ as the base metric for appropriate horizontal
data. While this is sometimes possible, it is noteworthy that the
`natural' choice of line bundle---$L$ pulled back over its own total
space, equipped with the induced Hermitian structure---does {\em
not}\/ fit into this framework. Indeed, if $(L,h)$ is non-flat, then
the curvature of the pullback bundle, computed with respect to an {\em
arbitrary}\/ bundle-adapted metric on the disk bundle, does not have
constant eigenvalues.
\end{remark}

For data as in Lemma~\ref{lemma:atom}, the functions $P$, $Q$, and $R$
are written explicitly as follows: Let the constant scalar curvature
of the Hodge manifold $(M^m,g)$ be~$\sigma$, and let $p:(L,h)\to(M,g)$
be the Hermitian line bundle with curvature
$\gamma=-2\pi\omega$. For each positive integer~$k$, the bundle
$(L^k,h^k)$ has curvature $k\gamma$. Fix $\alpha>0$, and equip $M$
with the metric $g_M=2\pi\alpha\,g$. Together with the compatible
interval $[0,\infty)$, these data are $\sigma$-constant, and
\begin{equation}
\label{eqn:PQR}
Q(\mum)=\left(1-{k\over\alpha} \mum\right)^m,\qquad
R(\mum)={{\alpha\sigma}\over{\alpha-k \mum}},\qquad
P(\mum)=2\sigma\left(1-{k\over\alpha} \mum\right)^{m-1}.
\end{equation}

Since every positive or negative Einstein-K\"ahler metric is
(homothetic to) a Hodge metric, Lemma~\ref{lemma:atom} encompasses
most prior examples of horizontal data to which Theorems
\ref{thm:A}~and \ref{thm:B} may be applied. It is reasonable to ask
whether or not a Hodge metric of constant scalar curvature is
``essentially'' Einstein-K\"ahler. The answer is ``no,'' even
discounting product metrics, homogeneous metrics, and the
like. LeBrun~\cite{LeBrun2} constructed scalar-flat metrics on certain
blown-up ruled surfaces; the K\"ahler classes of these metrics depend
on a real parameter, and the class is rational when the parameter
value is rational. Higher-dimensional examples include the manifolds
obtained from $\P^{2k+1}$ by blowing up a pair of skew $\P^k$'s. Such
a manifold is a $\P^1$-bundle over $\P^k\times\P^k$, and a K\"ahler
class is determined by the area of a $\P^1$~fibre and by the areas of
lines in each of the base factors. If the latter are equal, and are a
rational multiple of the area of a fibre, then the K\"ahler class is
proportional to a Hodge class, and by~\cite{Hwang} (Remark~5.3,
p.~584) is represented by a metric of constant scalar curvature.

The family $g_M(\mum)$ of horizontal metrics arising from
$\sigma$-constant horizontal data need not be homothetic. A simple way
to arrange this is to take a suitable product of data arising from
Lemma~\ref{lemma:atom}:

\begin{lemma}
\label{lemma:product}
Let $\{p_j:(L_j,h_j)\to(M_j,g_j),\ [0,\infty)\}$, $j=1,\ldots,n$, be
$\sigma$-constant horizontal data. Then the line bundle
$$p:L=\bigotimes_{j=1}^n \pi_j^*L_j\longrightarrow
M=M_1\times\cdots\times M_n,\qquad
\pi_j:M\to M_j\quad\hbox{the projection},$$ 
equipped with the induced metrics, is $\sigma$-constant.
\end{lemma}

The proof is immediate, and is left to the reader. The corresponding
function~$Q$ is a product of terms as in equation~(\ref{eqn:PQR}),
while $R$ is a sum of such terms. 

If $(M_j,g_j)$ are Einstein-K\"ahler, then the construction just
described gives rise to complete metrics on tensor products of
pluri-canonical and pluri-anticanonical bundles as in
Corollaries~\ref{cor:A1} and \ref{cor:A2}.

\subsubsection*{Remarks about compact Einstein-K\"ahler metrics}

By Yau and Aubin's solution of the Calabi conjecture, simply-connected
Einstein-K\"ahler manifolds of non-positive scalar curvature are
plentiful, the simplest examples being smooth complete intersection
varieties of degree $k\geq N+1$ and dimension $m\geq2$ in the complex
projective space~$\P^N$.

Work of Tian and Yau~\cite{TY} and Tian~\cite{Tian} shows that, with
precisely two exceptions, every compact complex {\em surface} with positive
first Chern class admits an Einstein-K\"ahler metric.

\begin{remark}
\label{remark:metrics}
(Deformation of positive Einstein-K\"ahler manifolds)\quad In higher
dimensions, work of Nadel~\cite{Nadel} and Siu~\cite{Siu} shows that
if $m/2\leq k\leq m+1$, then the Fermat hypersurface of degree~$k$ in
$\P^{m+1}$ admits an Einstein-K\"ahler metric with positive Ricci
curvature. Nadel also proved existence of positive Einstein-K\"ahler
metrics on certain branched coverings of projective space.  These
examples can be deformed; this is an immediate consequence of a
Theorem of LeBrun and Simanca~\cite{LBS1} regarding deformations of
extremal K\"ahler metrics. Precisely, if $(M,J,g_M)$ is a positive
Einstein-K\"ahler manifold with no non-trivial holomorphic vector
fields, their result implies that every sufficiently small deformation
of the complex manifold $(M,J)$ also admits an Einstein-K\"ahler
metric.
\end{remark}

If $(M,J)$ is a smooth hypersurface of degree $k\leq m+1$ in
$\P^{m+1}$, then every small deformation is realized by a hypersurface
in the same projective space. Consequently, there is a moduli space of
dimension ${{m+k+1}\choose k}-(m+2)^2+1$ (the dimension of the space
of monomials minus the dimension of $\aut\P^{m+1}$) consisting of
Einstein-K\"ahler structures (complex structure and compatible
K\"ahler metric) on the smooth manifold underlying the Fermat
hypersurface. Since for $m\geq3$ every smooth, irreducible
hypersurface in $\P^{m+1}$ has $h^{1,1}=b^2=1$ (by the Lefschetz
Hyperplane Theorem), every holomorphic line bundle over such a
manifold admits an Hermitian structure for which the line bundle data
are $\sigma$-constant.

\subsection{Vector Bundles}

This section describes some examples of data satisfying the hypotheses
of Theorem~\ref{thm:C}. As mentioned in the introduction,
$\sigma$-constancy for a vector bundle of rank $n>1$ is a strong
condition, nowhere nearly as flexible as the corresponding notion for
line bundles. This is expected, since the base metric~$g_D$ lives on a
space of dimension~$d$, while it is necessary to control $d+n-1$
curvature eigenvalues. The examples here are all chosen so that the
exceptional divisor of the blow-up is particularly simple, either a
homogeneous space, or else a product $D\times\P^{n-1}$. Nonetheless,
there are interesting metrics, many of which seem to be new.

\subsubsection*{Homogeneous vector bundles}

Let $D$ be a {\em compact}, homogeneous K\"ahlerian manifold, and fix a
maximal compact group $K\subset\aut^0(D)$, endowed with a bi-invariant
measure of unit volume. The group $K$ is unique up to conjugacy in 
$\aut^0(D)$. Furthermore, $K$ acts transitively on $D$, and
every de~Rham class (in particular, every K\"ahler class) contains a 
{\em unique} $K$-invariant representative, obtained from an arbitrary
representative by averaging. 

\begin{remark}
\label{remark:homog}
When $D$ is simply-connected, it is known that $D$ is rational, and
that every holomorphic vector bundle over $D$ is homogeneous.
Generally, a compact, homogeneous K\"ahler manifold is the K\"ahler
product of a flat torus and a rational homogeneous space, see
Matsushima~{\rm\cite{Matsushima2}}, or Besse~\cite{Besse} for a more
detailed expository treatment.
\end{remark}

The Ricci form of a $K$-invariant K\"ahler metric is $K$-invariant and
represents $2\pi c_1(D)$. In other words, there is only one $K$-invariant 
Ricci form, and its eigenvalues, with respect to a $K$-invariant K\"ahler
form are $K$-invariant functions, i.e.\ constants. 

Suppose $p:E\to D$ is a homogeneous holomorphic vector bundle, i.e.\
is induced by a representation of $K$ on $GL(n,\C)$, and that the
ruled manifold $M=\P(E)$ is $K$-homogeneous (which happens, e.g., when
$E$ is irreducible). Then the tautological bundle $p:L\to M$ is a
homogeneous line bundle, and by averaging over~$K$ it is clear that
for every K\"ahler class on~$M$, there exists a K\"ahler form
$\omega_M$ representing the chosen class, and an Hermitian structure
in~$L$, such that the data $p:(L,h)\to(M,g_M)$ are $\sigma$-constant.

A tubular neighbourhood of the zero section of~$L$ is obtained by
blowing up a tubular neighbourhood of the zero section of~$E$, so it
is clear that these neighbourhoods are either both pseudoconvex or
both not pseudoconvex. But in~$L$, pseudoconvexity is equivalent to
non-positivity of the curvature.  Thus, under the assumptions of
Corollary~\ref{cor:homog}, the hypotheses of Theorem~\ref{thm:C} are
satisfied, so the total space of~$E$ (or the disk subbundle) admits
complete metrics of constant scalar curvature.

Finally, consider the $K$-invariant Ricci form $\rho_M$
on~$M$. The restriction of $\rho_M$ to a fibre of $\pi:M\to D$
is the Ricci form of a Fubini-Study metric, so the form
$\rho_M-n\gamma$ vanishes on $\ker\pi_*$, hence is pulled back
from~$D$. Provided this form is non-positive, there is an
Einstein-K\"ahler metric on the total space of~$E$ or the disk
subbundle. We shall not pursue this avenue; a detailed (partial)
classification of cohomogeneity-one Einstein-K\"ahler metrics has
been accomplished by Dancer and Wang~\cite{DW}, see also Podesta and
Spiro~\cite{PS}. 

\subsubsection*{Sums of line bundles}

Let $E=\Lambda\otimes\C^n$ be a sum of $n$~copies of an Hermitian line
bundle over a base space such that $(\Lambda,h)\to(D,g_D)$ is
$\sigma$-constant. Then the projectivization is a product manifold
$M=D\times\P^{n-1}$, and the tautological bundle is tensor product
$\tau_E=\pi_1^*\Lambda\otimes\pi_2^*\O_{\P^{n-1}}(-1)$.
It is clear that the induced Hermitian structure on $L=\tau_E$ is
$\sigma$-constant with respect to the obvious metric on~$M$.

As a partial complement to this remark, observe that if $L_1$~and
$L_2$ are holomorphic line bundles over a {\em compact}\/
manifold~$D$, and if $L_1^*\otimes L_2$ admits a non-trivial
holomorphic section, then there do not exist metrics $g$~and $h$
such that the data $(L_1\oplus L_2,h)\to(D,g)$ are
$\sigma$-constant. This generalizes the remark about Hirzebruch
surfaces that was made following Theorem~\ref{thm:C}. As in that case,
the projective bundle $M=\P(L_1\oplus L_2)$ has non-reductive
automorphism group, hence (by a theorem of Lichnerowicz) cannot admit
a K\"ahler metric of constant scalar curvature. Considerations of this
type rather seriously restrict the possibility of finding sums of line
bundles to which Theorem~\ref{thm:C} applies. For example, if 
$D=\P^d$ (or more generally, is a compact, irreducible, rank-one
Hermitian symmetric space), then the {\em only}\/ sums of line bundles
satisfying the hypotheses of Theorem~\ref{thm:C} are as in the previous
paragraph. 

\subsubsection*{Stable bundles over curves}

Many of the observations below regarding ruled manifolds over curves
were made by Fujiki in the context of seeking extremal K\"ahler
metrics on ruled manifolds, see~\cite{Fujiki}.

Let $C$ be a smooth, compact Riemann surface of genus at least two,
and let $E \to C$ be a holomorphic vector bundle of rank~$n$ and
degree~$k$. We will use a theorem of Narasimhan-Seshadri~\cite{NS}, as
formulated by Atiyah-Bott~\cite{AB}, to show that when $E$ is stable
there is an Hermitian metric in~$E$ satisfying the hypotheses of
Theorem~\ref{thm:C}.

Equip $C$ with the (unique up to isometry) K\"ahler form~$\omega_C$ of
unit area and constant Gaussian curvature, so that $c_1(E) =
k\omega_C$.  In \cite[\S\S6 and 8]{AB} it is explained that stable
holomorphic bundles~$E$ correspond to irreducible representations
$\Gamma_R \to U(n)$, where $\Gamma_R$ is a central extension of the
fundamental group:
$$ 0 \to \R \to \Gamma_R \to \pi_1(C) \to 0.$$
It follows that if $E$ is stable then there is a family of metrics
with constant Ricci eigenvalues on $\P(E)$.  Indeed, because the
universal cover~$\Delta$ of~$C$ is Stein (and contractable) the
universal cover of $\P(E)$ may be identified holomorphically with
$\Delta \times \P(\C^n)$.  Equipping this with a product metric
$\hbox{Poincar\'e}\times\hbox{Fubini-Study}$, it follows that if $E$
is stable then $M=\P(E)$ is the quotient $\Delta \times \P(\C^n)/
\pi_1(C)$ where the action is by {\em isometries}. Scaling the two
factors gives a two-parameter family of K\"ahler forms on~$M$ whose
Ricci eigenvalues are constant, and whose eigenbundles are the
vertical tangent bundle of~$M$ and its orthogonal complement. Thus
$h^{1,1}(M)\geq2$.

The Leray-Hirsch Theorem (see, e.g.,~\cite{Kob2}, pp.~31ff.) implies
that, as a module over $H^*(C,\R)$, the cohomology ring $H^*(M,\R)$ is
generated by the first Chern class~$\zeta$ of the tautological
bundle~$\tau_E$ subject to the relation 
\begin{equation}
\label{eqn:degree}
0=\zeta^n-c_1(E)\,\zeta^{n-1}
=\zeta^{n}-k\omega_C\,\zeta^{n-1}.
\end{equation} 
In particular, $h^2(M)=2$, so by the observations made in the previous
paragraph, every two-dimensional cohomology class is represented by a
form whose pullback to $\Delta\times\P^{n-1}$ is a (possibly
indefinite) linear combination of the Poincar\'e form and the
Fubini-Study form. K\"ahler classes are exactly those classes whose
representatives pull back to positive combinations of these metrics.
Fix a K\"ahler class on~$M$, and let $g_M$ be the distinguished
representative. Every holomorphic line bundle $p:L\to(M,g_M)$ admits
an Hermitian structure~$h$, unique up to scaling, whose curvature form
is a combination of the Poincar\'e and Fubini-Study forms; thus the
data $p:(L,h)\to(M,g_M)$ are $\sigma$-constant.

\begin{lemma}
\label{lemma:taut-form}
Let $E\to C$ be a holomorphic vector bundle of rank~$n$ and degree~$k$
over a compact Riemann surface of genus $g\geq2$. Assume $C$~and
$M=\P(E)$ are equipped with metrics as above, and let $\omega_F$
denote the push-forward to $M$ of the integral Fubini-Study form on
$\Delta\times\P^{n-1}$. Then the curvature form $\gamma$ of the
tautological bundle of~$E$ and the Ricci form $\rho_M$ are given by
\begin{equation}
\label{eqn:taut-form}
{1\over{2\pi}}\gamma={k\over n}\omega_C-\omega_F,\qquad
{1\over{2\pi}}\rho_M=(2-2g)\,\omega_C+n\,\omega_F.
\end{equation}
\end{lemma}

\begin{proof}
Write $\gamma=k_1\,\omega_C+k_2\,\omega_F$. Then $k_2=-1$ since the
restriction of $\tau_E$ to a fibre is $\O_{\P^{n-1}}(-1)$.  By
equation~(\ref{eqn:degree}) and a short calculation, $k_1=k/n$.  To
see the Ricci form is as claimed, pull back to the universal cover,
where the metric is a product, and recall that $\omega_C$ has unit
area. 
\end{proof}

If~$L\to M$ is a line bundle whose first Chern class is non-positive,
then the hypotheses of Theorems \ref{thm:A}~and \ref{thm:B} are
satisfied. Further, by Lemma~\ref{lemma:taut-form} the tautological bundle
$L=\tau_E$ has non-positive first Chern class iff $k=\deg E\leq0$. In
this event, the hypotheses of Theorem~\ref{thm:C} are satisfied, and
the total space of~$E$ (or the disk subbundle) admits complete
K\"ahler metrics of constant scalar curvature. Finally,
Lemma~\ref{lemma:taut-form} implies
$\rho_M+n\gamma=2\pi(2g-2+k)\,\omega_C$, so if $k=2-2g$, then by
Theorem~\ref{thm:C} the total space of~$E$ admits a complete,
Ricci-flat K\"ahler metric, while if $k<2-2g$, the disk subbundle
of~$E$ supports a complete Einstein-K\"ahler metric of negative
curvature.

\subsection{Metrics of Finite Fibre Area}
\label{finitev}

Metrics of finite fibre area arise in the momentum construction when
the momentum interval is bounded. In this case the most convenient
normalization is to take the momentum interval~$I$ to be symmetric
about $\mum=0$, rather than insisting that $\inf I=0$. The requirement
that $\omega_M(\mum)=\omega_M-\mum\gamma$ be positive-definite for all
$\mum\in I$ no longer imposes a sign restriction on~$\gamma$. The main
difficulty is that the problem of constructing a complete,
finite-volume K\"ahler metric of constant scalar curvature is
overdetermined: There are two boundary conditions---namely
$\varphi=0$, and $|\varphi'|=0$~or $2$---that must be satisfied at
each end of the momentum interval, and only three parameters (the
initial conditions, and the value of the scalar
curvature). Consequently, it is to be expected that the K\"ahler class
of~$g_M$ will be restricted by existence of a complete metric with
constant scalar curvature. 

\subsubsection*{Metrics on compact manifolds}

Theorem~\ref{thm:D} below is an existence result for K\"ahler metrics
of constant (positive) scalar curvature on certain compact manifolds
to which the momentum construction is applicable. It was proven
in~\cite{Hwang}, and is included here for two purposes: To suggest the
type of theorem to be expected for non-compact (but finite-volume)
metrics, and to indicate the parts of the proof that generalize with
no extra effort. The important philosophical point is that on a
compact manifold, it is not generally the case that every K\"ahler
class is represented by a metric of constant scalar curvature.

\begin{thm}
\label{thm:D}
Let $(M,g_M)$ be a product of positive Einstein-K\"ahler manifolds,
each having $b_2=1$, and let $p:L\to M$ be a holomorphic line bundle
whose first Chern class is strictly indefinite. Then the completion
$\widehat{L}=\P(L\oplus{\bf1})$ admits a K\"ahler metric of constant
scalar curvature. In fact, the set of K\"ahler classes containing such
a metric is a real-algebraic family that separates the K\"ahler cone.
\end{thm}

\begin{proof}
(Sketch) The first step is to establish that under the hypotheses of
Theorem~\ref{thm:D}, every K\"ahler class on the compact
manifold~$\widehat{L}$ is represented by an extremal K\"ahler metric
(in the sense of Calabi). For present purposes, this may be taken to
mean that the scalar curvature is an affine function of the momentum
coordinate~$\mum$, i.e.\ that the gradient field of the scalar
curvature is a global holomorphic vector field. The terminology comes
from a variational problem having the latter property as its
Euler-Lagrange equation, see~\cite{Cal1}. A K\"ahler class containing
an extremal representative is an {\em extremal class}.

Let $Q$~and $R$ be defined as in Section~\ref{section:KS}, and let
$I=[-b,b]$ be the momentum interval.  Proceeding backward, set
$\sigma(\mum)=\sigma_0+\sigma_1\mum$ and attempt to solve the boundary
value problem 
\begin{equation}
\label{eqn:extr}
(\varphi Q)''(\mum)
=2Q(\mum)\Bigl(R(\mum)-\sigma_0-\sigma_1\mum\Bigr);
\qquad \varphi(\pm b)=0,\quad \varphi'(\pm b)=\mp2.
\end{equation}
The values of $\sigma_0$~and $\sigma_1$ are determined uniquely by
$Q$, $R$, and $b$ (Lemma~\ref{lemma:bd-values} below), and an
elementary (but slightly involved) root counting argument shows that
$(\varphi Q)''$ vanishes at most twice, so that $\varphi Q$---which is
positive near the endpoints of~$I$---is positive on $(-b,b)$. Thus
equation~(\ref{eqn:extr}) determines a momentum profile whose induced
metric is extremal, proving that every K\"ahler class on~$\widehat{L}$
contains an extremal representative (possibly with non-constant scalar
curvature). The root counting argument rests crucially on the fact
that the curvature of the base metric is non-negative.

The proof is completed by determining conditions under which
$\sigma_1=0$, see also~(\ref{eqn:f=0}) below. This is accomplished by
expressing $\sigma_1$ as a polynomial in the curvature and Ricci
eigenvalues of the horizontal data, so that the ``variables'' are
exactly the parameters controlling the K\"ahler class of the base
metric. Roughly, the top coefficient changes sign as the parameters
vary, so $\sigma_1$ changes sign on the K\"ahler cone, hence vanishes
on a real-algebraic hypersurface that separates the cone. The
corresponding metrics have constant scalar curvature.
\end{proof}

There are substantial difficulties in extending Theorem~\ref{thm:D} to
compact manifolds when the base curvature is not positive. Perhaps the
greatest, found by T{\o}nnesen-Friedman~\cite{CTF}, is that the set of
extremal classes is not obviously the entire K\"ahler cone. In
particular, on a ruled surface whose base has genus at least two, the
set of classes for which the momentum construction yields an extremal
metric is {\em not} generally the entire K\"ahler cone. This
potentially complicates the final portion of the argument, since
$\sigma_1$ may vanish for certain choices of eigenvalues, but the
relevant parameters may not correspond to extremal classes. 

\subsubsection*{Non-compact metrics of finite fibre area}

By analogy with the compact case, it is desirable to search among a
family of metrics whose scalar curvature may be non-constant. The
natural extension is to the class of {\em formally extremal} metrics,
by definition those whose scalar curvature is an affine function of
the moment map. As before, the hope is to find, for each pair of
boundary conditions, an affine function $\sigma_0+\sigma_1\mum$ such
that the function~$\varphi$ satisfying
\begin{equation}
\label{eqn:bd-values}
R(\mum)-{1\over{2Q(\mum)}}(\varphi Q)''(\mum)=\sigma_0+\sigma_1\mum
\end{equation}
matches the given boundary conditions and is non-negative. Matching
the boundary conditions is easy linear algebra:

\begin{lemma}
\label{lemma:bd-values}
Let $\{p:(L,h)\to(M,g_M),[-b,b]\}$ be $\sigma$-constant horizontal
data for some $b>0$. For each pair of boundary values $\varphi'(\pm
b)$, there exists a unique choice of $\sigma_0$ and $\sigma_1$ such
that the function $\varphi$ defined by {\rm(\ref{eqn:bd-values})} has
the given boundary derivatives and satisfies $\varphi(\pm b)=0$.
\end{lemma}

\begin{proof}
For $n\geq0$, set 
$$\a{n} =\int_{-b}^b x^n Q(x)\,dx,\qquad \A{n} 
= -\mum^n(\varphi Q)'(\mum)\Big|_{-b}^b+\int_{-b}^b x^n(QR)(x)\,dx.$$
Observe that $\a{n} /\a{0} $ is the $n$th moment (in the sense of
probability) of the moment map~$\mum$, computed with respect to the
symplectic measure, and that $\a{2} \a{0} -\a{1} ^2 >0$ by the
Cauchy-Schwarz Inequality. Integrating~(\ref{eqn:bd-values}), using
the given boundary conditions at $\mum=-b$, yields
\begin{eqnarray}
(\varphi Q)'(\mum) & = & (\varphi Q)'(-b)
+2\int_{-b}^\mum \Bigl(R(x)-\sigma_0-\sigma_1x\Bigr)Q(x)\,dx, 
\nonumber \\
&& \label{eqn:bd-values1} \\
(\varphi Q)(\mum) & = & (\varphi Q)'(-b)(\mum+b)
+2\int_{-b}^\mum(\mum-x)\Bigl(R(x)-\sigma_0-\sigma_1x\Bigr)Q(x)\,dx. 
\nonumber 
\end{eqnarray}
Setting $\mum=b$ and subtracting $b$~times the first equation from the
second leads to the system
$$\left[\matrix{\a{0} & \a{1} \cr \a{1} & \a{2} \cr}\right]
\left[\matrix{\sigma_0 \cr \sigma_1 \cr}\right]
=\left[\matrix{\A{0} \cr \A{1} \cr}\right].$$
This system has a unique solution since the coefficient matrix is
non-singular.
\end{proof}

This ``boundary-matching'' result is the only part of the proof of
Theorem~\ref{thm:D} that generalizes immediately, and the results of
T{\o}nnesen-Friedman suggest that there are genuine geometric
complications in less restricted settings.  It is appropriate to
remark here that 
\begin{equation}
\label{eqn:f=0}
\sigma_1=0\qquad {\rm iff}\qquad \a{0} \A{1} -\a{1} \A{0} =0.
\end{equation}
In the compact case, this is exactly the condition that the Futaki
character of a K\"ahler class vanishes; in the present situation, it
is to be expected that there is an invariant for ``finite volume
K\"ahler classes'' on certain non-compact manifolds. However, it is
not immediately obvious whether or not the profile is positive on a
neighbourhood of $\mum=\pm b$, regardless of whether or not
$\sigma_1=0$; in the compact case, this is automatic because the
boundary derivatives are non-zero.

\subsection{Other Constructions of Bundle-Adapted Metrics}
\label{section:compare}

Several authors have constructed bundle-adapted metrics from various
points of view. With the exception of LeBrun, who worked over curves,
the authors mentioned below have used the following curvature
hypotheses, either implicitly or explicitly:

\begin{description}
\item{(i)} The eigenvalues of the curvature endomorphism~$\Beta$ are
constant on~$M$;
\item{(ii)} The eigenvalues of the Ricci endomorphism~$\varrho$ are
constant on~$M$;
\item{(iii)} At each point of~$M$, $\Beta$~and $\varrho$ are
simultaneously diagonalizable.
\end{description}

For brevity, data satisfying these conditions are said to be {\em
$\rho$-constant}. Data which are $\rho$-constant are clearly
$\sigma$-constant. Roughly, the distinction is between assuming an
endomorphism has constant trace and assuming its eigenvalues are
constant. 

\subsubsection*{Historical survey}

The list below is in approximate chronological order,
but does not necessarily follow lines of development back to their
earliest discernible origins.

Calabi~\cite{Cal0} used distortion potential functions to construct
complete Einstein-K\"ahler metrics in line bundles over an
Einstein-K\"ahler base, and in the cotangent bundle of~$\P^d$, and
used the same method (in~\cite{Cal1}) to construct compact
``extremal'' K\"ahler metrics of non-constant scalar curvature.

Koiso and Sakane~\cite{KS1,KS2} used the momentum map as a coordinate
to construct compact Einstein-K\"ahler metrics of real cohomogeneity
one.  Their work followed Sakane~\cite{Sakane}, likely inspired by the
work of B\'erard-Bergery~\cite{Bergery} on compact, non-homogeneous
Einstein metrics. One salient point of their work is an explicit
interpretation of vanishing of the Futaki invariant as an
integrability condition for an ODE.  At about the same time,
Mabuchi~\cite{Mab} gave a more symplectic proof of existence of
Einstein-K\"ahler metrics on the same spaces considered by Koiso and
Sakane. Mabuchi also gave a satisfactory interpretation of the Futaki
invariant in terms of convex geometry.

A symplectic approach was also taken by LeBrun~\cite{Lebrun}, and his
formalism was generalized by Pedersen and Poon~\cite{PP}. LeBrun
assumed that the base $M$ was a curve, but allowed the profile to be a
general (positive) function on $I\times M$. (Correspondingly the
dependence upon $\mum$ of $\omega_M(\mum)$ is not necessarily affine,
and $d\omega_\varphi=0$ is an additional condition.)  He showed that
{\em every} $S^1$-invariant scalar-flat K\"ahler metric in complex
dimension 2 is locally described by a pair of functions $u$ and $w$ on
$I\times M$ satisfying certain partial differential equations. His $w$
is essentially the reciprocal of our $\varphi$.  What is remarkable is
that these equations are tractable when $\dim M=1$ and lead, for
example, to the construction of scalar-flat K\"ahler metrics on
(certain) blow-ups of ruled surfaces.

Pedersen and Poon allowed $\dim M >1$ and worked with torus bundles,
but made assumptions that reduce their very complicated system of
equations to the ODEs studied in this paper. (Indeed in most of their
examples, $M$ is Einstein--K\"ahler and $L$ is a (possibly fractional)
power of the canonical bundle.)  Their examples include metrics of
constant scalar curvature on line and disk bundles over projective
spaces.  In particular they showed that there exist K\"ahler metrics
of zero scalar curvature on the total space of $\O(-m) \to \P_n$
provided that $m \geq n$. In fact this restriction on $m$ is
unnecessary, as was shown by Simanca~\cite{Sim1}, also using the
Calabi ansatz.

The methods of Koiso and Sakane were used in~\cite{Hwang}, and
independently by Guan~\cite{Guan1}, to extend Calabi's families of
extremal metrics. Hwang (\cite{Hwang}, p.~564) attributed the
construction to Koiso and Sakane, overlooking the fact that Calabi
(\cite{Cal0}, p.~281, equation~(4.9), for example) had written the
moment map and distance function in the manner of
equation~(\ref{eqn:mu}) above. However, Calabi seems to have made the
observation in passing, and did not emphasize the use of momentum
coordinates.

Koiso~\cite{Koiso} and Guan~\cite{Guan2} used the method of Koiso and
Sakane to study Hamilton's Ricci flow for K\"ahler metrics and how it
may fail to converge; the concept of a {\em quasi-Einstein} K\"ahler
metric is introduced in the latter two papers. A quasi-Einstein
K\"ahler metric is the K\"ahlerian analogue of a {\em Ricci soliton},
introduced by R.~Hamilton.

Dancer and Wang~\cite{DW} and Podesta and Spiro~\cite{PS} have
independently used the techniques of Koiso and Sakane to obtain a
partial classification of Einstein-K\"ahler metrics having real
hypersurface orbits under the action of the isometry group.

T{\o}nnesen-Friedman~\cite{CTF} used the Calabi ansatz to study
existence of extremal K\"ahler metrics on some ruled surfaces. Some of
her examples are of great interest for the following reason. If
$(M,J)$ is a compact complex manifold that admits an extremal K\"ahler
metric in some K\"ahler class, then there are two properties that
might generally be hoped for:
\begin{itemize}
\item {\em Every}\/ K\"ahler class contains an extremal
representative;
\item Extremal metrics in a fixed class on~$M$ are unique up to the
action of the automorphism group.
\end{itemize}
T{\o}nnesen-Friedman found families of complex surfaces (admitting
extremal metrics) for which at least one of the preceeding statements
fails.

The papers of Calabi~\cite{Cal0}, Koiso-Sakane~\cite{KS1}, and
Pedersen-Poon~\cite{PP} develop versions of the momentum construction
from the points of view of K\"ahlerian, Riemannian, and symplectic
geometry, respectively. The remainder of this section gives short
summaries of each method (in the notation of the present paper where
applicable), including dictionaries between these constructions and
the momentum construction as presented here.

\subsubsection*{The Calabi construction: K\"ahler distortion potentials}

Let $p:(E,h)\to(D,g_D)$ be an Hermitian holomorphic vector bundle of
rank~$n$ over a K\"ahler manifold of dimension~$d$. In a coordinate
chart $V\subset D$ over which $E$~is trivial, there exist holomorphic
coordinates $z=(z^\alpha)$, $\alpha=1,\ldots,d$ for~$D$ and fibre
coordinates---i.e.\ local holomorphic sections of~$E$---denoted
$\zeta=(\zeta^i)$, $i=1,\ldots,n$, which together constitute a chart
for the total space of~$E$.  The Hermitian structure is given in this
chart by an Hermitian matrix-valued function
$H=(H_{i\bar\jmath}):V\to\C^{n\times n}$, and the norm squared
function~$r$ is given by $r=\zeta^i\bar{\zeta}^j\,H_{i\bar\jmath}(z)$.
The canonical Hermitian connection is given locally by the form
$\Theta=H^{-1}\del H$, and the curvature form is
$\Omega=\delbar\Theta$.

Let $E'\subset E$ be an invariant (disk) subbundle, in the sense of
Definition~\ref{def:invariant}.  For every strictly subharmonic
function $\widehat{F}:E'\to\R$, $\omega=p^*\omega_D+\ddbar
\widehat{F}$ is a K\"ahler form on $E'$. Following Calabi~\cite{Cal1},
$\widehat{F}$ is called a {\em distortion potential} for $\omega$. The
idea is to construct a K\"ahler form on~$E'$ by lifting the
``horizontally supported'' form $\omega_D$ and adding a form which is
non-degenerate in the fibre directions. Some care is required since
$\ddbar\widehat{F}$ has horizontal components, of course.

The Calabi ansatz is to choose $\widehat{F}=F(r)$. Thus
\begin{equation}
\label{eqn:form}
\omega=p^*\omega_D+\ddbar F(r).
\end{equation}
Explicitly, the assumption is that the level sets of $\widehat{F}$
coincide with the level sets of~$r$. A short calculation gives
$$\omega=\Bigl(F'(r)+rF''(r)\Bigr)
\left[\sqrt{-1}{{\del r}\over r}\wedge{{\delbar r}\over r}\right]
+p^*\omega_D-\Bigl(rF'(r)\Bigr)\,\gamma,$$ where the form~$\gamma$ is
the ``bi-Hermitian curvature form'' of~$(E,h)$, namely, the curvature
form of the tautological line bundle $L\to\P(E)$, pulled back to~$E$.

Calabi~\cite{Cal0} searches for Einstein-K\"ahler metrics, i.e.\
K\"ahler metrics with $\rho=k\omega$, whose (locally defined) K\"ahler
potential function~$\Phi$ satisfies the complex Monge-Amp\`ere
equation
$$\det(\Phi_{\alpha\bar\beta})=|{\rm hol}|^2\,e^{-k\Phi}$$ with $|{\rm
hol}|$ the absolute value of a non-vanishing (local) holomorphic
function on~$D$. Assuming further that~$g_D$ is Einstein-K\"ahler
and~$(E,h)$ is a line bundle with constant curvature, Calabi reduces
existence of an Einstein-K\"ahler metric to an ODE, and solves this
equation in terms of a polynomial function which, in our notation,
is~$\varphi Q$. Finally, he gives criteria (in terms of the horizontal
data) for completeness of the resulting metric, as in
Theorem~\ref{thm:A} above. Similar methods are applied to the
cotangent bundle of~$\P^d$.

Under Calabi's curvature hypotheses, the scalar curvature
of~(\ref{eqn:form}) is constant on the level sets of~$r$, and is
therefore specified by an ordinary differential expression in~$F$.
Unfortunately, this expression is fourth-order and fully nonlinear, so
it is not trivial (or always possible) to write the metrics explicitly
or to understand their geometry. In the case of a single fibre (i.e. a
bundle over a point), the scalar curvature of the K\"ahler
form~(\ref{eqn:form}) is 
$$\sigma=-e^{-\psi(r)}\Bigl(\psi'(r)+r\psi''(r)\Bigr),\qquad
\psi(r)=\log\Bigl(F'(r)+rF''(r)\Bigr);$$ 
generally, curvature terms from the bundle $(E,h)$ and base metric
$(D,g_D)$ enter. Solving the equation $\sigma=c$ for~$F$ is sometimes
possible, see Simanca~\cite{Sim1}, who treats line bundles over a
complex projective space.

\subsubsection*{The Koiso-Sakane ansatz}

The approach of Koiso and Sakane~\cite{KS1} is more in the spirit of
Riemannian geometry than K\"ahler geometry. Divide the
$\C^\times$-bundle $\lx$ by the natural circle action, obtaining
$(0,\infty)\times M$, then seek a function $s:(0,\infty)\to\R$ and a
family $\{g_s\}$ of Riemannian metrics on~$M$ such that forming the
warped product metric $ds^2+g_s$ on $(0,\infty)\times M$, lifting back
to~$\lx$, and taking the Hermitian metric ``$ds^2+(ds\circ J)^2+g_s$''
yields a K\"ahler metric. They determine that the family $g_s$ must be
of the form $g_M-\mum\,B$ (a special case of the Duistermaat-Heckman
theorem), and change coordinates so that the moment map $\mum$ is the
``independent variable.''

As we have emphasized, the geometry of $\omega$ is more easily
extracted from the momentum description than from the distortion
potential description due to the interplay between~$\varphi$
and~$\mum$.  To reiterate, the ``independent variable''~$t$ depends
only on the complex structure of the line bundle~$L$. By contrast, the
``independent variable'' $\mum=\mom(t)$ {\em depends on the choice of
profile}; the coordinate in which the metric is described is unknown
data at the outset. It is remarkable that the non-linearity of the
scalar curvature---which is locked into the description when using a
fixed holomorphic coordinate system---is absorbed into the unknown
momentum map, leaving only a ``universal'' second-order {\em linear}
operator. Holomorphic coordinates do not always provide the most
transparent geometric description, even for K\"ahler metrics.

\subsubsection*{The Hamiltonian constructions of LeBrun, Pedersen and Poon} 

Holomorphic coordinates are also suppressed in this approach.  Start
with a complex manifold $M$, an interval $I$ and a circle bundle
$L'\to I\times M$, equipped with a $U(1)$-connection. The aim is to
give~$L'$ an $S^1$-invariant Riemannian metric that is K\"ahler with
respect to a complex structure to be constructed from data on $I\times
M$. Let $\Theta$ be the connection $1$-form on~$L'$, $\w$~a positive
function on $I\times M$, and $\mum$~an affine parameter on~$I$. The
holomorphic structure of~$L'$ is specified by taking $\w\,d\mum +
i\Theta$, along with the $(1,0)$ forms pulled back from~$M$, to be of
type $(1,0)$. The ansatz for the K\"ahler form is
$$
\omega = d\mum\wedge\Theta + \omega',
$$
where $\omega'$ is a family of $(1,1)$-forms on $M$, parametrized by
$I$. It follows from these definitions that $1/\w$ is the
length-squared of $X$, the  generator of the $U(1)$-action on $L'$.
In general the conditions of integrability and closure lead to a
complicated system of equations; these can be found in \cite{PP},
where the generalization to torus bundles was also considered.

Things are simpler when $\dim M =1$, the case considered by LeBrun
\cite{Lebrun}. He wrote 
$z=x+iy$ for a complex coordinate on $M$ and put
$$\omega = d\mum\wedge \Theta +
\w e^u dx\wedge dy$$
for some smooth $u$ on $I\times M$. The conditions that $\omega$ be
closed and the almost-complex structure integrable are given by the
PDE
\begin{equation}
\label{eqn:lebrun1}
\w_{xx} + \w_{yy} + (\w e^u)_{\mum\mum} =0.
\end{equation}
The condition $\sigma=0$ is equivalent to the so-called $SU(\infty)$
Toda lattice equation 
\begin{equation}
\label{eqn:lebrun2}
u_{xx} + u_{yy} +(e^u)_{\mum\mum} =0.
\end{equation}
The metric arises from the Calabi ansatz when $\w$ is independent of
$x$ and $y$, since $\w^{-1}=\varphi(\mum)$. In this event,
equation~(\ref{eqn:lebrun1}) forces $\w e^u$ to be an affine function
of $\mum$, precisely in accord with the momentum construction.

When $\omega$ is highly symmetrical it can happen that there is a
choice of circle action that yields a description in terms of the
momentum construction, while another choice of circle action does
not. For example, the scalar-flat `Burns metric'~$\omega$ on
$\widetilde{\C}^2$, the blow-up of $\C^2$ at the origin, arises in the
momentum construction from the identification
$\widetilde{\C}^2=\O(-1) \to \P^1$. However, the Burns metric is
$U(2)$-invariant, and may be symplectically reduced with respect to
the $S^1$-action given in standard coordinates on $\C^2$ by $(Z_1,Z_2)
\mapsto (e^{i\theta}Z_1,Z_2)$.  LeBrun calculated the functions
$\w$~and $u$, and found $e^u=2\mum$ but a complicated expression for
$\w$, not independent of $(x,y)$. In particular, with this choice of
horizontal data $\omega$ does {\em not} arise from the momentum
construction!

\end{document}